\documentclass[11pt]{amsart}

\usepackage[english]{babel}

\newtheorem{theorem}{Theorem}[section]
\newtheorem{lemma}[theorem]{Lemma}
\newtheorem{corollary}[theorem]{Corollary}
\newtheorem{proposition}[theorem]{Proposition}
\newtheorem{conjecture}{Conjecture}
\newtheorem{claim}[theorem]{Claim}
\newtheorem{problem}{Problem}

\theoremstyle{remark}
\newtheorem{rem}[theorem]{Remark}
\theoremstyle{definition}
\newtheorem{definition}[theorem]{Definition}

\def\bR{\mathbb R} \def\bQ{\mathbb Q} \def\bN{\mathbb N} \def\bZ{\mathbb Z} \def\bP{\mathbb P}
\def\brp{\mathbb R \mathbb P} \def\bcp{\mathbb C \mathbb P}
\def\mC{\mathcal{C}} \def\mA{\mathcal{A}} \def\mP{\mathcal{P}} \def\mL{\mathcal{L}} \def\mS{\mathcal{S}} \def\mH{\mathcal{H}} \def\mD{\mathcal{D}} \def\mE{\mathcal{E}} \def\mZ{\mathcal{Z}}
\def\bfE{\mbox{{\bf X}}}\def\bfS{\mbox{{\bf S}}}\def\bfR{\mbox{{\bf R}}}\def\bfB{\mbox{{\bf E}}}\def\bfH{\mbox{{\bf H}}}\def\bfG{\mbox{{\bf G}}}\def\bfr{\mbox{{\bf r}}}\def\bfK{\mbox{{\bf K}}}\def\bfC{\mbox{{\bf C}}}\def\bfN{\mbox{{\bf N}}}

\DeclareMathOperator{\im}{Im}
\DeclareMathOperator{\del}{\partial}
\def\cratio{\bfr}
\def\ra{\rightarrow}
\def\co{\colon\thinspace}
\def\Vlabel{\mbox{normal}}

\newcommand{\defin}[1]{\textit{#1}}

\begin{document}

\title{Toric geometry of convex quadrilaterals}
\author{Eveline Legendre}
\date{}
\thanks{This paper is part of the author's Ph.D. thesis. The author would like to thank her supervisors, Vestislav Apostolov and Paul Gauduchon, for hours of enlightening discussions. She would also like to thank David Calderbank for useful suggestions. }

\address{D\'epartement de Math\'ematiques, UQAM \\
C.P. 8888, Succ. Centre-ville Montr\'eal (Qu\'ebec), H3C 3P8, Canada\\
Centre de Math\'ematiques Laurent Schwartz, \'Ecole Polytechnique, 91128 Palaiseau, France}
\email{eveline.legendre@cirget.ca}

\keywords{extremal K\"ahler metrics, toric $4$--orbifolds, Hamiltonian $2$--forms}
\subjclass{Primary 53C25; Secondary 58E11}

\begin{abstract}
We provide an explicit resolution of the Abreu equation on convex labeled quadrilaterals. This confirms a conjecture of Donaldson in this particular case and implies a complete classification of the explicit toric K\"ahler--Einstein and toric Sasaki--Einstein metrics constructed in \cite{H2FII,exMS,other2}. As a byproduct, we obtain a wealth of extremal toric (complex) orbi-surfaces, including K\"ahler--Einstein ones, and show that for a toric orbi-surface with 4 fixed points of the torus action, the vanishing of the Futaki invariant is a necessary and sufficient condition for the existence of K\"ahler metric with constant scalar curvature. Our results also provide explicit examples of relative $K$--unstable toric orbi-surfaces that do not admit extremal metrics.
\end{abstract}

\maketitle

\section{Introduction}\label{secFUTintro}

This paper classifies toric K\"ahler metrics admitting Hamiltonian $2$--forms on $4$--dimensional toric orbifolds. Apostolov, Calderbank and Gauduchon~\cite{H2FI} pointed out that these $2$--forms underpin known explicit constructions of extremal K\"ahler metrics (in the sense of Calabi), see e.g. Calabi~\cite{calabi} and Bryant~\cite{bryant}, which situate our work into the more central problem of finding (explicit) extremal metrics on compact symplectic toric manifolds and orbifolds. A closely related problem in Sasakian geometry is the study of compatible toric Sasaki metrics which are transversally extremal~\cite{BGS}. The common feature of these two problems is that, by using the toric assumption, they can be reduced to a quasi-linear $4$--th order PDE on a (convex, compact, simple) polytope in $\bR^n$, which we solve explicitly for convex quadrilaterals in $\bR^2$.\\

We now describe the more general setting of the problem, following the work of Guillemin~\cite{guillMET}, Abreu~\cite{abreu} and Donaldson~\cite{don:scalar}. Let $\Delta\subset \bR^n$ be a convex, compact, simple polytope and $u= \{u_1, \dots, u_d\}$ a set of vectors in $(\bR^n)^*$ inward to $\Delta$ and respectively normal to the facets $F_1$, ... $F_d$ of $\Delta$. By slight abuse of notation, we shall refer to $(\Delta,u)$ as a \defin{labeled polytope}. In the case when there is a lattice $\Lambda\subset \bR^n$ with $u_i\in\Lambda$, there is a unique positive integer $m_i\in\bN$ such that $\frac{1}{m_i}u_i$ is a primitive element of $\Lambda$. Then $(\Delta, m_1,\dots, m_d)$ is a \defin{rational labeled polytope} in the sense of Lerman--Tolman~\cite{LT:orbiToric} and describes a compact symplectic toric orbifold.\\

 The set $\mS(\Delta, u)$ of \defin{symplectic potentials} consists of smooth strictly convex functions, $G\in C^{\infty}(\mathring{\Delta})$, whose inverse Hessian $$\bfH=(H_{ij}) = \left(\frac{\del^2 G}{\del \mu_i \del \mu_j}\right)^{-1}$$ is smooth on $\Delta$, positive definite on the interior of any face and satisfies, for every $y$ in the interior of the facet $F_i\subset \Delta$,
\begin{equation}\label{condCOMPACTIFonH}
  \bfH_y (u_i, \cdot) =0\;\;\;\mbox{ and }\;\;\; d\bfH_y (u_i, u_i) =2u_i.
\end{equation} These expressions use the standard identifications $(\bR^n)^* \simeq \bR^n$, $T_{x}\bR^n \simeq \bR^n$. The boundary conditions (\ref{condCOMPACTIFonH}) are elaborated in~\cite{H2FII} and are equivalent to the boundary conditions in~\cite{don:estimate} and~\cite{abreuOrbifold}. When $(\Delta,u)$ is a rational labeled polytope corresponding to a symplectic toric orbifold, $\mS(\Delta,u)$ parameterizes the space of compatible K\"ahler metrics. As computed by Abreu~\cite{abreu}, the scalar curvature of the metric associated to a symplectic potential $G\in \mS(\Delta,u)$ is the pull-back by the moment map of the function
\begin{equation}\label{abreuForm}
S(G)=-\sum_{ij} \frac{\del^2 H_{ij}}{\del \mu_i \del \mu_j}.
\end{equation}

The \defin{extremal affine function} $\zeta_{(\Delta,u)}$ is the $L^2$--projection (with respect to the euclidian measure) of $S(G)$ to the finite dimensional space of affine-linear functions on $\Delta$. It turns out that $\zeta_{(\Delta,u)}$ is independent of the symplectic potential $G\in S(\Delta,u)$ and may also be defined as the solution of a linear system depending only on $(\Delta,u)$, see \S~\ref{sectSURVEYmetric}. If $(\Delta,u)$ is a rational labeled polytope associated to a symplectic toric orbifold, then the symplectic gradient of the pull-back by the moment map of $\zeta_{(\Delta,u)}$ is the extremal vector field~\cite{futakimabuchi}.\\

The general problem we are interested in is then
\begin{problem}\label{problem1} Given a labeled polytope $(\Delta,u)$, is there a symplectic potential $G\in \mS(\Delta,u)$ satisfying the extremal K\"ahler equation
\begin{equation}\label{extrem} S(G) = -\sum_{ij} \frac{\del^2 H_{ij}}{\del \mu_i \del \mu_j}= \zeta_{(\Delta,u)} \;?
\end{equation} If so, can one find it explicitly ?\end{problem}

Bryant~\cite{bryant} and Abreu~\cite{abreuOrbifold} showed that on labeled simplices, the solution of this problem is given (up to an additive affine-linear function) by \begin{equation}\label{stadSOLSimplexe}
G= \frac{1}{2}\left(\sum_{i=1}^d \ell_i \log \ell_i - \ell_{\infty} \log \ell_{\infty}\right)\end{equation} where $\ell_i(\cdot)=\langle \cdot ,u_i \rangle-\lambda_i$ is an affine-linear function such that $F_i \subset \ell_i^{-1}(0)$ and $\ell_{\infty} = \sum_{i=1}^{d} \ell_i$. Explicit solutions are also known to exist for trapezoids corresponding to toric Hirzebruch surfaces, see \cite{calabi,abreu}.

Motivated by the conjectured link \cite{yau,tian:conj,don:scalar} between the existence problem for K\"ahler metrics of constant scalar curvature lying in an integer K\"ahler class and stability of the Kodaira embedding of the corresponding polarized variety, Donaldson~\cite{don:scalar} gave a precise conjecture for the existence part of Problem \ref{problem1}. It is expressed in terms of positivity of a linear functional $\mL_{\Delta,u}$, called the {\it relative Futaki functional} in this paper, over a suitable space of convex functions (on $\Delta$). On a labeled polytope $(\Delta,u)$, the relative Futaki functional is defined by
\begin{equation} \label{RelatFutakiFonct}
 \mL_{\Delta,u} (f) =\int_{\del\Delta}fd\nu - \frac{1}{2}\int_{\Delta}f\zeta_{(\Delta,u)}dv,
\end{equation}
where $dv$ is an euclidian measure on $\Delta$ and $d\nu$ is a measure on any facet $F_i$ defined by $u_i\wedge d\nu= -dv$.

\begin{definition} A labeled polytope $(\Delta, u)$ is \defin{analytically relatively $K$--stable} {\it with respect to toric degenerations} if the associated relative Futaki functional $\mL_{\Delta,u}$ is non-negative on any convex continuous piecewise affine-linear function on $\Delta$, and vanishes if and only if the function is affine-linear.
\end{definition}

\begin{conjecture}\label{conjecturDoSz} Let $(\Delta, u)$ be a labeled polytope. There is a solution $G\in\mS(\Delta, u)$ to the equation~(\ref{extrem}) if and only if $(\Delta, u)$ is analytically relatively $K$--stable with respect to toric degenerations.
\end{conjecture}

Donaldson proved Conjecture~\ref{conjecturDoSz} for polygons ($n=2$) when $\zeta_{(\Delta,u)}$ is constant, by using the continuity method~\cite{don:scalar,don:estimate,don:extMcond,don:csc}. Zhou--Zhu proved that the existence of a solution of~(\ref{extrem}) implies analytical relative $K$--stability, see~\cite[Proposition 2.2]{zz}\footnote{The result is stated for toric manifolds but the proof goes through for any labeled polytope.}.\\

In this paper, we consider the case where $(\Delta, u)$ is a labeled convex quadrilateral ($\zeta_{(\Delta,u)}$ is not required to be constant). We show that there is a compatible K\"ahler metric admitting a Hamiltonian $2$--form on any symplectic toric orbifolds whose moment polytope is a quadrilateral. Using the work of~\cite{wsd,H2FI} to separate variables in equation~(\ref{extrem}) in this case, we provide an explicit solution by means of elementary techniques in the case when $\zeta_{(\Delta,u)}$ is {\it equipoised} on $\Delta$.

\begin{definition} \label{invOrthoFct}
Let $\Delta$ be a quadrilateral with vertices $s_1,\dots,s_4$, such that $s_1$ is not consecutive to $s_3$. We say that the affine function $f$ is \defin{equipoised} on $\Delta$ if $$\sum_{i=1}^4 (-1)^i f(s_i) =0.$$
\end{definition}
More precisely, our main result is
\begin{theorem}~\label{mainTHEOparall}
Let $(\Delta, u)$ be a labeled convex quadrilateral with equipoised extremal affine function $\zeta_{(\Delta,u)}$. Then there exist two polynomials of degree at most $4$, $A(x)$ and $B(y)$, from which one can construct explicitly a $S^2\bR^2$--valued function $\bfH_{A,B}=(H_{ij})$ satisfying~(\ref{condCOMPACTIFonH}) and~(\ref{extrem}).

Moreover, $(\Delta, u)$ is analytically relatively $K$--stable with respect to toric degenerations if and only if $\bfH_{A,B}$ is the inverse Hessian of a solution $G_{A,B}\in\mS(\Delta, u)$ of~(\ref{extrem}) which happens if and only if $\bfH_{A,B}$ is positive definite\footnote{This condition is expressed as $A(x)$ and $B(y)$ being positive on certain intervals.}. In particular, if~(\ref{extrem}) admits a solution in $\mS(\Delta, u)$ this solution is given by $G_{A,B}$.  \end{theorem}

Note that our result provides a computable condition of relative $K$--stability for the labeled polytopes we consider and gives an explicit solution for the constant scalar curvature equation on labeled quadrilaterals. Indeed, constant functions are equipoised on any quadrilateral. Another natural class of examples is given by toric weakly Bocher-flat K\"ahler orbi-surfaces see~\cite{wsd} which must be either K\"ahler--Einstein, Bochner flat (given by~(\ref{stadSOLSimplexe}) on weighted projective planes, see~\cite{bryant}) or given by the explicit solution of the above form on certain labeled quadrilaterals with equipoised extremal function, see Remark~\ref{remWBF}. Nevertheless, there exist labeled convex quadrilaterals whose extremal affine function is {\it not} equipoised and which admit a solution to Problem~\ref{problem1}. Indeed, it follows from~\cite{don:extMcond} that the set of inward normals to a convex quadrilateral $\Delta$ for which $(\Delta,u)$ admits a solution to Problem~\ref{problem1} is open in $(\bR^4)^*$ and, when $\Delta$ is not a parallelogram, intersects the hyper-surface of normals $u$ for which $(\Delta,u)$ has an equipoised extremal affine function (see Theorem~\ref{propCONEsol} below).\\

One might wonder how restrictive it is to require the extremal affine function of a labeled quadrilateral to be equipoised. A key point to answer this (as well as to prove our main result) is the observation that fixing the polytope $\Delta$ and varying the inward normals $u$ the coefficients of the polynomials $A(x)$, $B(y)$ depend {\it linearly} on $u$. More precisely, for a convex quadrilateral $\Delta$, denote by $\bfN(\Delta)$ the $4$--dimensional cone of inward normals $u=(u_1,u_2, u_3, u_4)$ associated to the facets of $\Delta$ and its subset $\bfB(\Delta)$ (respectively $\bfC(\Delta)$) defined by the condition that $\zeta_{(\Delta,u)}$ is equipoised (respectively constant). Define $\bfB^+(\Delta) $ as the set of normals $u\in \bfB(\Delta)$ such that $(\Delta,u)$ is relatively analytically $K$--stable with respect to toric degenerations. Studying these sets leads to the following result.

\begin{theorem} \label{propCONEsol} Let $\Delta$ be a convex labeled quadrilateral which is not a parallelogram. $\bfC(\Delta)$ is a codimension-one sub-cone of $\bfB(\Delta)$ which is itself a codimension-one sub-cone of $\bfN(\Delta)$. Furthermore,  $\bfB^+(\Delta)$ contains $\bfC(\Delta)$ and is a non-empty open subset of $\bfB(\Delta)$, which is proper if $\Delta$ is not a trapezoid. Finally, there is a $1$--dimensional cone $\bfK(\Delta)\subset \bfC(\Delta)$ such that the corresponding solutions of $(3)$ define (homothetic) K\"ahler--Einstein metrics on $\mathring{\Delta} \times \bR^2$. \end{theorem}

The above theorem implies that for convex labeled quadrilaterals, the condition that $\zeta_{(\Delta,u)}$ is constant implies the $K$-stability of $(\Delta,u)$. It also provides an effective parametrization of the K\"ahler--Einstein solutions found in~\cite{exMS,other2,H2FII} in terms of classes of affine-equivalent convex quadrilaterals.\\

One geometric application concerns symplectic toric $4$--orbifolds. Using the well-known correspondence~\cite{delzant:corres,LT:orbiToric} between compact symplectic toric orbifolds and (convex, compact, simple) rational labeled polytopes, we obtain as a corollary of Theorems~\ref{mainTHEOparall} and~\ref{propCONEsol} that a symplectic toric $4$--orbifold whose labeled polytope is a quadrilateral with equipoised extremal affine function admits (an explicit) compatible extremal metric if and only if its labeled polytope is analytically relatively $K$--stable with respect to toric degenerations. Moreover, it admits a constant scalar curvature metric if and only if the Futaki invariant vanishes (i.e. $\zeta_{\Delta,u}$ is constant).

The question whether or not a polytope $\Delta$ is \defin{of rational type}, that is, corresponds to the image of the moment map of a compact symplectic toric orbifold, preludes the problem of constructing examples systematically. We answer this question for convex polygons by giving an explicit criterion in terms of the cross-ratio of the elements of $\brp^1$ corresponding to the facets of $\Delta$, see Theorem~\ref{crossratioTHEO}. We then show that for a {\it strongly rational polytope} $\Delta$ (i.e whose vertices lie in a lattice) the linear constraints on the normals introducing the cones $\bfB(\Delta)$, $\bfC(\Delta)$ and $\bfK(\Delta)$ of Theorem~\ref{propCONEsol} have rational coefficients, thus obtaining a wealth of examples of both extremal and unstable toric orbifolds. As corollary we obtain the following existence result.

\begin{corollary}\label{theoExistKE} Let $\Delta$ be a convex quadrilateral. If $\Delta$ is strongly rational then, up to homothetic transformations and finite orbifold coverings, there exists a unique compact K\"ahler--Einstein toric orbifold having $\Delta$ as moment polytope.\end{corollary}

Another geometric application concerns Sasaki toric $5$--manifolds. There is a correspondence between connected compact co-oriented contact toric manifolds of Reeb type and strictly convex good polyhedral cones~\cite{L:contactCONVEX,L:contactToric,contactNOTE}. Moreover, the Reeb vector field $X$ defines a labeled polytope $(\Delta_{X},u_{X})$ which is not rational unless the Reeb vector field is quasi-regular, see~\cite{contactNOTE} and~\S\ref{sectSasakSurvey}. Theorem~\ref{mainTHEOparall} leads to the following.

\begin{corollary}\label{corodimension5}Let $(N^5,g,T^3)$ be a compact Sasaki toric 5-manifold with constant scalar curvature, Reeb vector field $X$ and momentum cone with $4$ facets. Then, the corresponding transversal labeled polytope $(\Delta_X, u_X)$ is a K-stable quadrilateral with respect to toric degenerations and $\zeta_{\Delta,u}$ constant. In particular, the metric $g$ is explicitly determined in terms of two polynomials of degree at most $3$.
\end{corollary}

Note that explicit Sasaki--Einstein metrics of this type have been found in~\cite{exMS,other2} and toric Sasaki--Einstein metrics have been systematically studied in~\cite{futakiChoOno,futakiOnoWang}.\\

In order to prove Theorem~\ref{mainTHEOparall}, we consider the cases of \defin{generic} labeled quadrilaterals (that is, quadrilaterals which have no parallel edges) and labeled trapezoids.

In the first case, we use the notion of orthotoric K\"ahler structures introduced by Apostolov--Calderbank--Gauduchon in~\cite{wsd}. We refer to Section~\ref{orthoDEFNandPARAM} for definition and properties of such structures. Loosely speaking, these metrics provide a separation of variables in the expression of the scalar curvature. On compact $4$--orbifolds, they also come with a Hamiltonian toric action with three or four fixed points. The case of three fixed points corresponds to triangles, that is, weighted projective spaces for which orthotoric K\"ahler metrics have been studied in great detail~\cite{bryant,abreuOrbifold,H2FII}. We prove in \S\ref{sectGEN=ortho} that there exists an orthotoric K\"ahler metric compatible with the symplectic form of any symplectic toric orbifold associated to a generic rational labeled quadrilateral. Moreover, if such symplectic toric orbifold has a equipoised extremal affine function, then any extremal compatible K\"ahler metric must be orthotoric.

To cover the case of labeled trapezoids, we first notice that for labeled parallelograms Problem~\ref{problem1} is trivially solved by taking the product of solutions on labeled intervals given by~(\ref{stadSOLSimplexe}). In the orbifold case, these solutions correspond to product of extremal metrics on weighted projective lines. Notice also that any affine-linear function on a parallelogram is equipoised. We then introduce the notion of a \defin{Calabi toric metric} on a toric orbifold, naturally extending the construction of Calabi~\cite{calabi}. We show that any symplectic toric $4$--orbifold whose labeled moment polytope is a labeled trapezoid (but not a rectangle) admits a Calabi toric metric. Moreover, if such symplectic toric orbifold has a equipoised extremal affine function, we show that the extremal compatible K\"ahler metric (should it exist) must be Calabi toric.\\

As a byproduct, we obtain a classification of compact K\"ahler toric $4$--orbifolds admitting non-trivial Hamiltonian $2$--forms. These forms are defined and studied in~\cite{wsd,H2FI,H2FII}. For instance, it is shown in these works that orthotoric metrics are exactly those admitting Hamiltonian $2$--forms of maximal order and examples of non-trivial Hamiltonian $2$--forms on weighted projective spaces are given. We obtain the following classification: The moment polytope of a K\"ahler toric $4$--orbifold $(M,\omega, J, g,T)$ admitting a non-trivial Hamiltonian $2$--form is either a triangle (so that $(M,J)$ is a weighted projective plane) or a quadrilateral. In this latter case, the order of a non-trivial Hamiltonian $2$-forms is $2-p$ where $p$ is the number of pair of parallel edges. Moreover, $g$ is orthotoric if and only if $p=0$, $g$ is Calabi toric if and only if $p=1$ and $g$ is a product of metrics if and only if $p=2$.

\section{Toric orbifolds and extremal metrics -- a quick survey. }\label{sectSURVEY}
\subsection{Toric orbifolds and labeled polytopes.}\label{sectSURVEYorbifVSpolyt}\label{sectAFFPolyt}\label{assLabePolyt}

 A \defin{polytope} $\Delta$ in an affine space is a bounded set given as the intersection of a finite number of (closed) affine half-spaces, $\mH_i$ ($1\leq i\leq d$). We suppose that $d$ is minimal. The polytopes we study in this paper lie in the underlying affine space of $\mathfrak{t}^*$, the dual of the Lie algebra $\mathfrak{t}$ of a $n$--torus. Thus, for each $\mH_i$, we can specify a normal inward vector $u_i\in \mathfrak{t}$ and a point $x_i\in\mH_i$, so that $\mH_i= \{x \;|\;\langle x-x_i, u_i \rangle\geq 0\}$. The polytope $\Delta$ is then described via the \defin{defining equations}:
\begin{equation}\label{defEquatPolyt}
   \Delta =\{ x\in \mathfrak{t}^* \;|\;\langle x, u_i\rangle \geq \lambda_i\;\; \mbox{ for } i= 1,\dots, d \}\;\mbox{ with } \; \lambda_i= \langle x_i, u_i\rangle.
\end{equation}

We suppose that the \defin{interior} of the polytope, $\mathring{\Delta}$, is a non-empty open subset of $\mathfrak{t}^*$. A \defin{face} of $\Delta$ is a non-empty subset $F$, for which there exists $I\subset \{1,\dots,d\}$ such that $F=F_I= \cap_{i\in I}\mH_i\cap \Delta$. In particular, $\Delta=F_{\emptyset}$ is a face and all faces are close. The \defin{facets} of $\Delta$ are $n-1$--dimensional faces denoted $F_i=F_{\{i\}}$ while the \defin{vertices} are faces of $\Delta$ consisting of only one point. A polytope is \defin{simple} if each vertex is the intersection of $n$ distinct facets, where $n$ is the dimension of $\mathring{\Delta}$.

\textit{From now on, we assume all polytopes to be simple, compact and convex.}

\begin{definition}\label{defnRatpolyt}
A polytope $\Delta$ is \defin{rational with respect to a lattice $\Lambda\subset \mathfrak{t}$}, if each facet admits a normal vector lying in $\Lambda$, so we say that the pair $(\Delta,\Lambda)$ is rational. A polytope $\Delta$ is \defin{of rational type} if there exists a lattice $\Lambda \subset \mathfrak{t}^*$ such that the pair $(\Delta,\Lambda)$ is rational; it is \defin{strongly rational} if, up to translation, its vertices lie in a lattice.\end{definition}

\begin{definition}\label{RationLabelPolyt}
A \defin{labeled polytope} $(\Delta, u_1, \dots, u_d)$ is a polytope $\Delta\subset \mathfrak{t}^*$ endowed with an inward normal vector attached to each facet. We call these vectors \defin{$\Vlabel$s} and say that $F_i$ is \defin{labeled} by $u_i$. A \defin{rational labeled polytope}, $(\Delta, \Lambda, u_1, \dots, u_d)$, is a labeled polytope together with a lattice $\Lambda$ such that $u_i \in \Lambda$. Since the polytope is simple, the $\Vlabel$s of a rational labeled polytope span a sublattice $\Lambda_{\rm{min}}=\mbox{span}_{\bZ}\{u_1,\dots, u_d\}\subset \Lambda$.
\end{definition}

Recall that a vector $v\in \mathfrak{t}$ is \defin{primitive} with respect to a lattice $\Lambda$ if it generates $\bR v \cap\Lambda$. A rational pair $(\Delta,\Lambda)$ may canonically be viewed as a rational labeled polytope by using the primitive vectors as $\Vlabel$s.

\begin{rem} \label{notationTolman} The notion of rational labeled polytope was introduced by Lerman and Tolman in~\cite{LT:orbiToric}. The rational labeled polytope $(\Delta, \Lambda, m_1,\dots,m_d)$ is defined as a rational pair $(\Delta,\Lambda)$, together with a positive integer $m_i$ attached to each facet $F_i$. From our viewpoint, $(\Delta, \Lambda, m_1,\dots,m_d)$ corresponds to $(\Delta, \Lambda,m_1w_1,\dots,m_dw_d)$ where $w_i$ is the unique inward vector normal to $F_i$ which is primitive with respect to $\Lambda$.
\end{rem}

\begin{definition}\label{polDELZANT} A \defin{Delzant} polytope is a rational pair $(\Delta,\Lambda)$, such that any vertex of $\Delta$ is the intersection of the facets whose (primitive) $\Vlabel$s form a $\bZ$--basis of $\Lambda$.
\end{definition}

\begin{definition}\label{defnPolytEquiv}The polytopes $\Delta\subset \mathfrak{t}^*$ and $\Delta' \subset \mathfrak{s}^*$ are \defin{equivalent} if there exists an affine isomorphism $\phi$ between the respective affine spaces so that $\phi(\Delta)=\Delta'$. Two labeled polytopes are \defin{equivalent} if the underlying polytopes are equivalent via an affine map, $\phi\co \mathfrak{t}^*\ra \mathfrak{s}^*$, and whose differential adjoint $(d\phi)^*\co \mathfrak{s}\ra \mathfrak{t}$ exchange the $\Vlabel$s. Two rational labeled polytopes are \defin{equivalent} if they are equivalent as labeled polytopes via an affine map $\phi$ whose differential's adjoint $(d\phi)^*$ exchanges the respective lattices.\end{definition}

We use the definition of orbifolds appearing in~\cite{LT:orbiToric}. We only consider diffeomorphisms between orbifolds. In particular, these maps are {\it good} in the sense of~\cite[Chapter 4]{BG:book}. We refer to the latter reference for a detailed exposition about orbifolds and to~\cite{LT:orbiToric} for an exposition of results about Lie group actions on orbifolds.

Let $(M,\omega)$ be a compact symplectic orbifold and $T$ a torus with Lie algebra $\mathfrak{t}$. Denote by $\mathfrak{t}^*$ the dual vector space of $\mathfrak{t}$. A \defin{Hamiltonian action} of $T$ on $(M,\omega)$ is a faithful representation $\rho \co T \ra \mathrm{Symp}(M,\omega)$ together with a $T$--equivariant smooth map $\mu \co M \ra \mathfrak{t}^*$ satisfying $d\mu(\xi) = -\iota_{d\rho(\xi)}\omega$. If the dimension of $T$ is half the dimension of $M$, $(M,\omega,\mu,T,\rho)$ is a \defin{symplectic toric orbifold}. Two such orbifolds, $(M,\omega,\mu,T,\rho)$ and $(M',\omega',\mu',T',\rho')$, are \defin{equivalent} if there exists a symplectomorphism (of orbifolds) $\psi \co (M,\omega) \ra (M',\omega')$ and an isomorphism $h\co T\ra T'$ such that $\psi \circ \rho(t)\circ \psi^{-1} = \rho'(h(t)) $ for all $t \in T$. We shall omit the representation $\rho$ from the notation when no confusion is possible.

It is well-known, see \cite{atiyah,convexMoment,delzant:corres}, that the image of the moment map of a toric manifold is a convex polytope in $\mathfrak{t}^*$. This polytope is rational with respect to the lattice $\Lambda =\ker (\mathrm{exp}\co \mathfrak{t}\ra T)$ and satisfies the Delzant condition of Definition~\ref{polDELZANT}. In the case of orbifolds, Lerman--Tolman~\cite{LT:orbiToric} showed that $(\im \mu, \Lambda)$ is a simple rational polytope which is Delzant if and only if $M$ is non-singular and for every $p\in M$, the orbifold structure group of $p$, say $\Gamma_p$, only depends on the smallest face $F$ containing $\mu(p)$. More precisely, if $F$ is $\Delta$ itself (that is, if $\mu(p)$ lies in the interior of $\Delta$) then $\Gamma_p$ is trivial and if $F$ is a facet then $\Gamma_p$ is isomorphic to $\bZ/m_{F}\bZ$ for some integer $m_F$ called the \defin{label} of $F$. Thus, any toric orbifold $(M,\omega,\mu,T,\rho)$ naturally defines a labeled polytope, $(\im \mu, \Lambda, \{m_{F}\})$, in the sense of~\cite{LT:orbiToric}. This corresponds to a rational labeled polytope, $(\im \mu,\Lambda, \{u_{F}\})$, by taking, for each facet $F\subset \Delta$, the $\Vlabel$ $u_{F}=m_F w_F$ where $w_F$ is the primitive inward normal vector to $F$.

The Delzant--Lerman--Tolman correspondence,~\cite{delzant:corres,LT:orbiToric}, states that the symplectic toric orbifold is determined by its associated rational labeled polytope, up to a $T$--invariant symplectomorphism (of orbifolds). Conversely, any rational labeled polytope can be obtained from a symplectic toric orbifold, via an explicit construction called \defin{Delzant's construction}. Two Hamiltonian actions $(\mu,T,\rho)$, $(\mu',T',\rho')$ on a symplectic orbifold $(M,\omega)$ are equivalent if and only if the associated labeled polytopes are equivalent, see Definition~\ref{defnPolytEquiv}. In~\cite{kkp}, this statement is proved in the smooth case, with $T'=T$. The orbifold counterpart is formally the same with, in addition, special attention paid to $\Vlabel$s on one side and weights of the torus action on the other.\\

Two rational labeled polytopes $(\Delta, \Lambda',u)$, $(\Delta, \Lambda,u)$ with $\Lambda'\subset\Lambda$ corresponds (via the Delzant construction) to a finite orbifold covering with deck transformation group $\Lambda/\Lambda'$. More precisely, $\Lambda/\Lambda'\subset T'=\mathfrak{t}/\Lambda'$ is finite and acts by symplectomorphisms on the symplectic toric orbifold $(M',\omega',T')$ associated to $(\Delta, \Lambda',u)$. The quotient of $M'$ by $\Lambda/\Lambda'$ is then a symplectic toric orbifold with respect to the torus $T=\mathfrak{t}/\Lambda$ and is associated to $(\Delta, \Lambda,u)$ via the Delzant--Lerman--Tolman correspondence. For studying $T$--invariant K\"ahler metric it is then not restrictive to consider the minimal lattice $\Lambda_{\mathrm{min}}=\mbox{span}_{\bZ}\{u_1,\dots, u_d\}$. Minimal lattices corresponds to \defin{simply connected orbifolds}, see~\cite{H2FII} and~\cite{Painchaud}. We shall omit the lattice from the notation of a rational labeled polytope when using the minimal lattice.\\

The cohomology class of the symplectic form of a symplectic toric orbifold is rational if and only the associated polytope is strongly rational with respect to the lattice of circle subgroups of the torus. Indeed, this property is a corollary of the original construction of Delzant~\cite{delzant:corres} and holds for orbifolds, see~\cite{LT:orbiToric}.

\subsection{Compatible K\"ahler metrics and the extremal equation}\label{sectSURVEYmetric}\label{sectSURVEYmetricCPT}\label{sectSURVEYmetricEXT}

We consider a K\"ahler toric orbifold $(M^{2n},\omega, J, g, T, \mu)$, where $g$ is $T$-invariant, $\omega$-compatible K\"ahler metric and $J$ is a complex structure such that $g(J\cdot,\cdot)=\omega(\cdot,\cdot)$. We denote by $(\Delta,u_1,\dots, u_d)$ the associated rational labeled polytope. Recall,\cite{delzant:corres,LT:orbiToric}, that $\mathring{M}=\mu^{-1}(\mathring{\Delta})$ is the subset of $M$ where the torus acts freely. The K\"ahler metric provides a horizontal distribution for the principal $T$--bundle $\mu \co \mathring{M} \ra \mathring{\Delta}$ which is spanned by the vector fields $JX_u$, $u\in \mathfrak{t}=\mbox{Lie } T$. This gives an identification between the tangent space at any point of $\mathring{M}$ and $\mathfrak{t} \oplus \mathfrak{t}^*$. Usually, one chooses a basis $(e_1,\dots,e_n)$ of $\mathfrak{t}$ to identify $\mathring{M} \simeq \mathring{\Delta}\times T$ using the flows of the induced vector fields $X_{e_1}$, ... $X_{e_n}$, $JX_{e_1}$, ... $JX_{e_n}$ (which commutes thanks to the integrability of $J$). The action-angle coordinates on $\mathring{M}$ are local coordinates $(\mu_1,\dots,\mu_d,t_1,\dots,t_d)$ on $\mathring{M}$ such that $\mu_i = \langle \mu,e_i\rangle$ and $X_{e_i} =\frac{\del}{\del t_i}$. The differentials $dt_i$ are real-valued closed $1$--forms globally defined on $\mathring{M}$ as dual of $X_{e_i}$ (i.e $dt_i(X_{e_i}) =\delta_{ij}$ and $dt_i(JX_{e_j})=0$).

In the action-angle coordinates $(\mu_1,\dots,\mu_d,t_1,\dots,t_d)$, the symplectic form becomes $\omega =\sum_{i=1}^n d\mu_i\wedge dt_i.$ It is well-known \cite{guillMET} that a K\"ahler toric metric may be expressed with respect to these coordinates as follows:
\begin{align}\label{ActionAnglemetric}
g= \sum_{s,r} G_{rs}d\mu_r\otimes d\mu_s + H_{rs}dt_r\otimes dt_s,
\end{align}
where the matrix valued functions $(G_{rs})$ and $(H_{rs})$ are smooth on $\mathring{\Delta}$, symmetric, positive definite and inverse to each other. In particular, $g_\mathrm{red}= \sum G_{rs}d\mu_r\otimes d\mu_s$ is a Riemannian metric on $\mathring{\Delta}$. It may be more convenient to view these objects through the identification between tangent spaces of $\mathring{M}$ and $\mathfrak{t} \oplus \mathfrak{t}^*$ as above. Indeed, following~\cite{H2FII} we define the $S^2\mathfrak{t}^*$--valued function $\bfH \co \mathring{\Delta} \ra \mathfrak{t}^*\otimes\mathfrak{t}^*$ by $\bfH_{\mu(p)}(u,v)= g_p(X_u,X_v)$, and put $H_{rs}= \bfH(e_r,e_s)$. Similarly, $\bfG \co\mathring{\Delta} \ra \mathfrak{t}\otimes\mathfrak{t}$ is the metric $g_\mathrm{red}$ via the usual identification $T_{\nu}\mathring{\Delta}\simeq \mathfrak{t}^*$. Given the expression~(\ref{ActionAnglemetric}), the integrability of the complex structure $J$ is equivalent to the relation
\begin{equation}\label{integCond}
  \frac{\del}{\del\mu_j} G_{rs} = \frac{\del}{\del\mu_r} G_{js}
\end{equation}
or, equivalently, to the fact that $(G_{rs})$ is the Hessian of a potential $G \in C^{\infty}(\mathring{\Delta})$. \\

Necessary and sufficient conditions for a $S^2\mathfrak{t}^*$--valued function $\bfH$ to be induced by a globally defined K\"ahler toric metric on $M$ are established in \cite{abreuOrbifold,don:estimate,H2FII}. We will use in this paper the boundary conditions of \cite{H2FII}, which we recall below. For a face $F = F_I= \cap_{i\in I} F_i$ of $\Delta$, denote $\mathfrak{t}_F = \mbox{span}_{\bR}\{ u_i\,|\; i\in I\}$. Its annihilator in $\mathfrak{t}^*$, denoted $\mathfrak{t}_F^o$, is naturally identified with $(\mathfrak{t}/\mathfrak{t}_F)^*$.

\begin{proposition}\label{prop1H2FII}\cite[Proposition 1]{H2FII}
Let $\bfH$ be a positive definite $S^2\mathfrak{t}^*$--valued function on $\mathring{\Delta}$, whose inverse satisfies (\ref{integCond}). $\bfH$ comes from a K\"ahler metric on $M$ if and only if $\bfH$ is the restriction to $\mathring{\Delta}$ of a smooth $S^2\mathfrak{t}^*$--valued function on $\Delta$, still denoted by $\bfH$, which verifies the boundary conditions (\ref{condCOMPACTIFonH}) and such that the restriction of $\bfH$ to the interior of any face $F \subset \Delta$ is a positive definite $S^2(\mathfrak{t}/\mathfrak{t}_F)^*$--valued function.
\end{proposition}
Recall that the set of symplectic potentials $\mS(\Delta,u)$, defined in the introduction, is the space of smooth strictly convex functions on $\mathring{\Delta}$ for which $\bfH=(\mbox{Hess }G)^{-1}$ satisfies the conditions of Proposition~\ref{prop1H2FII}. Notice that the compactification condition (\ref{condCOMPACTIFonH}) uses the $\Vlabel$s and not the lattice. This agrees with the fact that different lattices (but identical $\Vlabel$s) lead to orbifolds admitting a common finite orbifold covering.

Abreu~\cite{abreu} computed the curvature of a compatible K\"ahler toric metric, $g$, in terms of its potential $G_g\in\mS(\Delta,u)$. More precisely, choosing a basis of $\mathfrak{t}=\mbox{ Lie } T$ as above and using the action-angle coordinates to express the metric as~(\ref{ActionAnglemetric}), the scalar curvature on $\mathring{M}$ is the pull-back by $\mu$ of the function $S(G_g)$, defined by~(\ref{abreuForm}) via the symplectic potential $G_g$ of $g$.

It is well-known, see e.g \cite{abreu}, that the metric $g$ is extremal if and only if $S(G_g)$ is an affine linear function on $\Delta$. In this case, this function must be equal to the \defin{extremal affine function} $\zeta_{(\Delta,u)}$ which we now define. Choosing a basis $(e_1,\dots, e_n)$ of $\mathfrak{t}$ gives a basis $\mu_0= 1$, $\mu_1=\langle e_1,\cdot\rangle$, $\dots$, $\mu_n=\langle e_n,\cdot\rangle $ of affine-linear functions. We define $\zeta_{(\Delta,u)}=\sum_{i=0}^n \zeta_i \mu_i$ where the vector $\zeta = (\zeta_0,\dots,\zeta_n)\in \bR^{n+1}$ is the unique solution of the linear system
\begin{equation}\label{systCHPext}
\begin{split}
\sum_{j=0}^n W_{ij}\, \zeta_j &= Z_i,\;\;\; i= 0,\dots, n\\
\mbox{with }\quad W_{ij} = \int_{\Delta} \mu_i\mu_jdv &\quad \mbox{ and }  \quad Z_i =2\!\int_{\del \Delta} \mu_i d\nu,
 \end{split}
\end{equation} where the volume form $dv=d\mu_1\wedge \dots \wedge d\mu_n$ and the measure $d\nu$ on $\del \Delta$ defined by the equality $u_j\wedge d\nu=-dv$ on the facet $F_j$. In other words, $\zeta_{(\Delta,u)}$ is determined by requiring that the linear functional (\ref{RelatFutakiFonct}) annihilates any affine-linear function $f$.

Equivalently, the extremal affine function is the $L^2(\Delta, dv)$--projection of $S(G_g)$, for any compatible K\"ahler toric metric $g$ to the finite dimensional space of affine-linear functions. Indeed, integrating~(\ref{abreuForm}) by part and using condition~(\ref{condCOMPACTIFonH}) we get
$$Z_i=\int_{\Delta} S(G_g) \mu_i dv =2\!\int_{\del \Delta} \mu_i d\nu.$$

\begin{rem} \label{dependenceExtchp2} Let $\bfH=(H_{ij})$ be any $S^2\mathfrak{t}^*$--valued function on $\Delta$ satisfying the compactification condition~(\ref{condCOMPACTIFonH}) but which is not necessarily positive definite. Using integration by parts as above, one can show that the $L^2(\Delta, dv)$--projection of the function $$S(\bfH)=-\sum_{i,j}\frac{\del^2 H_{ij}}{\del \mu_i \del \mu_j}$$ on the space of affine-linear functions on $\Delta$ is still equal to $\zeta_{(\Delta, u)}.
$\end{rem}

\subsubsection{Uniqueness}\label{uniqueness}

Guan~\cite{guan} showed the uniqueness, up to automorphisms, of compatible extremal K\"ahler toric metrics on a smooth compact symplectic toric manifold. In fact, he proved that a geodesic in the space of compatible $T$--invariant K\"ahler metrics corresponds, via the Moser Lemma, to a straight line in the space of symplectic potentials $\mS(\Delta, u)$. Then, any two K\"ahler toric metrics may be linked together by a geodesic and, by a well-known argument using the convexity of the Guan--Mabuchi--Simanca \defin{relative $K$--energy} $\mE$ over geodesics, any two extremal toric metrics must coincide up to automorphisms.

 Guan's proof can be recasted in terms of symplectic potentials following the work of~\cite{don:scalar}. Indeed, it is showed that the relative $K$--energy of a metric associated to a potential $G\in \mS(\Delta, u)$ is
$$\mE(G)=2\mL_{(\Delta,u)}(G) - \int_{\Delta} (\log\det \mbox{Hess}(G))dv$$
where $\mL_{(\Delta,u)}$ is the relative Futaki functional~(\ref{RelatFutakiFonct}). One can check that the critical points of $\mE$ are exactly the solutions of~(\ref{extrem}) by computing that \begin{equation}
  \begin{split}
  d\mE_G(f)&=2\mL_{(\Delta,u)}(f) - \int_{\Delta}\langle \bfH , \mbox{Hess} f\rangle dv\\
  &= -\int_{\Delta} \bigg(\zeta_{(\Delta,u)} +\sum_{i,j} \frac{\del^2H_{ij}}{\del \mu_i \mu_j}\bigg)f dv,
  \end{split}
\end{equation} where $\bfH=(\mbox{Hess} G)^{-1}$ and $\langle\cdot ,\cdot \rangle$ is the usual inner product of symmetric matrices
(i.e the trace of their product). Moreover, for a segment $G_t = tG_1 +(1-t) G_0$ in $\mS(\Delta, u)$, we compute
\begin{equation}\label{eq:2thdiffMabuchiEnergy}
  \frac{d^2}{dt^2}\mE(G_t)= \int_{\Delta} \langle \bfH_t (\bfG_1-\bfG_0),\bfH_t (\bfG_1-\bfG_0)\rangle dv,
\end{equation} where $\bfG_t$ denotes the Hessian of $G_t$ and $\bfH_t=\bfG_t^{-1}$. This implies that $\mE$ is convex along segments in $\mS(\Delta, u)$.
Hence, if $G_0$ and $G_1$ are two solutions of~(\ref{extrem}) then $\mE$ is constant along $G_t$ and $\frac{d^2}{dt^2}\mE(G_t)=0$. Due to~(\ref{eq:2thdiffMabuchiEnergy}),
this implies that $\bfG_1-\bfG_0=0$, that is, $G_1-G_0$ is affine-linear.

\subsection{Sasaki toric geometry.}\label{sectSasakSurvey}
\subsubsection{Contact toric manifolds.}

Recall that there is a correspondence between co-oriented compact connected contact manifolds and symplectic cones over compact manifolds. Indeed, to such a contact manifold, $(N^{2n+1}, \mD)$, one can naturally associate the symplectic cone $(\mD_+^o,\hat{\omega},\varsigma)$ where $\mD_+^o$ is a connected component of the annihilator in $T^*N$ of the contact distribution $\mD$ without the zero section, $\hat{\omega}=d\lambda$ is the restriction of the differential of the canonical Liouville form $\lambda$ of $T^*N$, and $\varsigma$ is the Liouville vector field defined as $\varsigma_{(p,\alpha)}= \frac{d}{ds}_{|_{s=0}} e^{s}\alpha_p$, so that $\mL_{\varsigma} \hat{\omega}= \hat{\omega}$.

\begin{definition}
A (compact) \defin{contact toric} manifold $(N^{2n+1}, \mD, \hat{T}^{n+1})$ is a co-oriented compact connected contact manifold $(N^{2n+1}, \mD)$ endowed with an effective action of a (maximal) torus $\hat{T}\hookrightarrow \mbox{Diff}(N)$ preserving the contact distribution $\mD$ and its co-orientation. Equivalently, the symplectic cone $(\mD_+^o,\hat{\omega},\varsigma)$ is toric with respect to the action of $\hat{T}$ and the Liouville vector field $\varsigma$ commutes with $\hat{T}$. We denote by $$\hat{\mu}\co \mD_+^o \ra \hat{\mathfrak{t}}^*=(\mbox{Lie} \hat{T})^*$$ the \defin{contact moment map}, which is the unique moment map of $(\mD_+^o,\hat{\omega}, \hat{T})$ which is homogeneous of order $1$ with respect to $\varsigma$, see~\cite{L:contactToric}.
\end{definition}

\begin{definition}
A polyhedral cone is \defin{good} with respect to a lattice $\Lambda$, if any facet $F_i$ has a normal vector lying in $\Lambda$ and, for any face $F_I= \cap_{i\in I} F_i$,
\begin{equation}\label{goodness}
\mbox{span}_{\bZ}\{\hat{u}_i\,|\, i\in I\}= \Lambda\cap \mbox{span}_{\bR}\{\hat{u}_i\,|\, i\in I\}
\end{equation}
where $\hat{u}_i$ denotes the normal vector to $F_i$ which is primitive in $\Lambda$.
\end{definition}

Lerman established \cite{L:contactCONVEX,L:contactToric} a correspondence between contact toric manifolds and good polyhedral cones. The image, $\im  \hat{\mu}$, of the contact moment map of a compact contact toric manifold $(N, \mD, \hat{T})$ does not contain $0$ and $\mC= \im \hat{\mu} \cup \{0\}$ is a convex, polyhedral cone which is good with respect to the lattice of circle subgroups, $\Lambda \subset \hat{\mathfrak{t}}$. $\mC$ is called \defin{the moment cone}. Conversely, any convex good polyhedral cone is the moment cone of a contact toric manifold, unique up to contactomorphisms.

\subsubsection{Reeb vector field and transversal K\"ahler toric geometry}

\begin{definition}
A contact toric manifold $(N, \mD,\hat{T})$ is of \defin{Reeb type} if there exists a contact form whose Reeb vector field is induced by an element of the Lie algebra of $\hat{T}$. Equivalently, $(N, \mD,\hat{T})$ is of Reeb type if the moment cone $\mC$ is strictly convex, that is, $\mC^*_+ = \{a\in \hat{\mathfrak{t}}\,|\; \forall x\in \mC\backslash \{0\}, \; \langle a,x\rangle>0\}\neq \{\emptyset\}$.
\end{definition}
Let $(N,\mD,\hat{T})$ be a contact toric manifold of Reeb type and with contact moment map $\hat{\mu}$. Let $d$ be the number of facets of $\mC$ and denote by $\hat{u}_1, \dots, \hat{u}_d$ the set of primitive vectors in $\Lambda$ labeling $\mC$. For any $b\in \mC^*_+$, there is a contact form, $\eta_b$, for which $X_b$ is a Reeb vector field. We denote by $P_b$ the hyperplane $P_b =\{ x\in \hat{\mathfrak{t}}^*\,|\, \langle b,x\rangle=1\}$; $\Delta_b$ the polytope $\Delta_b = \mC\cap P_b$ and $\varrho_b$ the quotient map $\varrho_b \co\hat{\mathfrak{t}}\ra \hat{\mathfrak{t}}/\bR b$. The polytope $\Delta_b$ is $n$--dimensional, simple and compact~\cite{contactNOTE}.

Moreover, if $b\in\mC^*_+$ and $\bR b\cap \Lambda \neq \{0\}$, it generates a circle subgroup of $\hat{T}$: $T_b = \bR b / (\bR b \cap \Lambda)$. This group acts on the cone $(\mD_+^o,\hat{\omega})$ via the inclusion $\iota_b \co T_b \hookrightarrow \hat{T}$, with moment map $\mu_b \co \mD_+^o \ra \bR$ given by $\iota_b^*\circ\hat{\mu}.$ The space of leaves $\mZ_b$ of the Reeb vector field $X_b$, or equivalently the symplectic reduction $\hat{\mu}_b^{-1}(1)/T_b$, is an orbifold naturally endowed with a symplectic form $\omega$ and a Hamiltonian torus action of $\hat{T}/T_b$ such that the associated rational labeled polytope is \begin{equation}\label{transversPolyt}
 (\Delta_b, \varrho_b(\hat{u}_1),\dots , \varrho_b(\hat{u}_d))
\end{equation}
with the affine identification $(\hat{\mathfrak{t}}/(\bR b))^*\simeq P_b$, see~\cite{contactNOTE}.

\begin{definition}
 A \defin{Sasaki toric} manifold $(N, \mD,g,\hat{T})$ is a Sasaki manifold whose underlying contact structure is toric with respect to $\hat{T}$ and whose metric is $\hat{T}$--invariant.
\end{definition}
The Sasaki toric manifold $(N, \mD,g,\hat{T})$ corresponds to the K\"ahler toric cone $(\mD^o_+,\hat{\omega},\hat{g},\varsigma, \hat{T},\hat{\mu})$ where the metric $\hat{g}$ is the cone metric of $g$, that is $\hat{g}$ is homogenous of order one with respect to $\varsigma$ and restricts to $g$ on the level set $\hat{g}(\varsigma,\varsigma)=1$. Recall that $\hat{g}$ is K\"ahler, toric and homogeneous of order $1$ with respect to the (holomorphic) Liouville vector field $\varsigma$. Notice that $J\varsigma$ is induced by an element $b\in \mC^*_+\subset\hat{\mathfrak{t}}$, so that $X_b= J\varsigma$ restricts to a Reeb vector field on $N\subset M$, where $N$ is seen as the subset of $M$ where $\hat{g}(\varsigma,\varsigma)=1$.

The transversal geometry of $g$ refers to the metric $\check{g}$, induced on $\mD$, that is, $g= \eta_b\otimes \eta_b +\check{g}$. Recall that $d\eta_b$ restricts to a symplectic structure on $\mD$, so that $(\mD,d\eta_b, \check{g})$ defines a transversal K\"ahler toric structure. When the space of leaves $\mZ_b$ of the Reeb vector field $X_b$ is an orbifold (that is, when $\bR b\cap \Lambda \neq\{0\}$), we can identify its tangent $V$--bundle (see~\cite{BG:book}) with $\mD$ so that $(\mZ_b,d\eta_b)$ is the symplectic toric orbifold associated to the labeled polytope~(\ref{transversPolyt}). We want to describe the transversal geometry of $g$ using this polytope, which is defined even when $\bR b\cap \Lambda \neq \{0\}$.

The K\"ahler toric structure $(\hat{\omega},\hat{g},J)$ is expressed on the set $\mathring{M}= \hat{\mu}^{-1}(\mathring{\mC}\backslash \{0\})$ where the torus acts freely with a $S^2\hat{\mathfrak{t}}$--valued function $\hat{\bfG}$, having an inverse $\hat{\bfH}$. The $S^2\hat{\mathfrak{t}}^*$--valued function $\hat{\bfH}$ must satisfy the Proposition~\ref{prop1H2FII}. Indeed, these conditions are local and correspond to the smooth extension of the metric over the singular orbits of the action, see the proof of~\cite[Proposition 1]{H2FII}.

\begin{lemma}\cite{L:contactToric} \label{lemReebH}
The functions $\hat{G}_{ij}$ and $\hat{H}_{ij}$ are homogeneous of respective orders $-1$ and $1$ with respect to $\varsigma$. Moreover, if $b\in \hat{\mathfrak{t}}$ induces the Reeb vector field $X_b=J\varsigma$, then for all $\hat{\mu}\in \mC$, $\hat{\bfH}_{\hat{\mu}}(b,\cdot)=\hat{\mu}$.
 \end{lemma} Recall that $\hat{\bfH}$ is a map $\hat{\bfH} :\mC \ra S^2\hat{\mathfrak{t}}^*$ and we denote by $\hat{\bfH}_{\hat{\mu}}$ its value at $\hat{\mu}\in \mC$.
\begin{proposition}\label{propHbSasak}
Let $\mC$ be a good cone with inward normals $\hat{u}_1$, ... $\hat{u}_d$. Let $b\in \mC^*_+$ and $\hat{\bfH}$ be a positive definite $S^2\hat{\mathfrak{t}}^*$--valued function, satisfying the conditions of Proposition~\ref{prop1H2FII} on the faces of $\mC$. For every $\hat{\mu}\in\mC$, put
$$\bfH^b_{\hat{\mu}} = \hat{\bfH}_{\hat{\mu}}- \frac{\hat{\mu} \otimes\hat{\mu} }{\langle \hat{\mu},b\rangle}.$$
Then $\bfH^b$ is a positive definite $S^2 (\hat{\mathfrak{t}}/\bR b)^*$--valued function, satisfying the conditions of Proposition~\ref{prop1H2FII} with respect to $(\Delta_b, \varrho_b(\hat{u}_1),\dots , \varrho_b(\hat{u}_d))$. Moreover, if $\bR b\cap \Lambda \neq \{0\}$ and $\hat{\bfH}$ is associated to the cone metric of a Sasaki metric, $g$, on $N$, then the transversal metric $\check{g}$ induced by $g$ on the orbifold $\mZ_b$ is associated to the restriction of $\bfH^b$ to $\Delta_b$.
\end{proposition}

\begin{proof} We denote by $[a]$, the equivalence class of $a$ in $\hat{\mathfrak{t}}/(\bR b)$. Due to Lemma~\ref{lemReebH},
$$\bfH^b_{\hat{\mu}} = \hat{\bfH}_{\hat{\mu}}- \frac{\hat{\bfH}_{\hat{\mu}}(b,\cdot) \otimes \hat{\bfH}_{\hat{\mu}}(b,\cdot) }{\hat{\bfH}_{\hat{\mu}}(b,b)}.$$ This $S^2\hat{\mathfrak{t}}^*$--valued function is well-defined on the quotient $\hat{\mathfrak{t}}/\bR b$ since $\bfH^b(b,\cdot) =0$. A facet $F_k$ lies in the annihilator of the attached normal vector $\hat{u}_k$ since its closure in $\mathfrak{t}^*$ contains $0$. Thus, for $y\in F_k$, using $\hat{\bfH}_y(\hat{u}_k,\cdot)=0$, we get $\bfH^b_y([u_k],\cdot) = \bfH^b_y(u_k,\cdot) = 0$ and
$$d_y\bfH^b([\hat{u}_k],[\hat{u}_k]) = d_y\hat{\bfH}(\hat{u}_k,\hat{u}_k) - 2\frac{\langle y,\hat{u}_k\rangle}{\langle y,b\rangle} \hat{u}_k  + \frac{\langle y,\hat{u}_k\rangle^2}{\langle y,b\rangle^2}b = 2\hat{u}_k.$$
The $S^2\hat{\mathfrak{t}}^*$--valued function associated to a toric metric $\check{g}$ is defined by
$$\bfH^{\check{g}}(a,c) = \check{g}(X_a,X_c)$$
for two vectors $a$, $c\in\hat{\mathfrak{t}}$. On the other hand, the K\"ahler metric $\check{g}$ induced on $\mZ_b$ by $\hat{g}$ is defined by $\check{g}_{[p]}([X],[X'])= \hat{g}_p( Z,Z')$ where $Z$ (resp. $Z'$) is the horizontal projection of $X$ (resp. $X'$) at $p$ and $[\;\!\cdot\;\!]$ denotes the equivalent class of points or vectors with respect to the local action generated by $X_b$. For $a\in \hat{\mathfrak{t}}$, the horizontal projection of $X_a$ is $$Z_{a}=X_a-\frac{\hat{\bfH}(a,b)}{\hat{\bfH}(b,b)}X_b.$$
 Then $\bfH^{\check{g}}([a],[c]) = \check{g}(X_{[a]},X_{[c]})=\check{g}([X_a],[X_c]) = \hat{g}_p(Z_a, Z_c) = \bfH^b([a],[c])$,
which concludes the proof.
\end{proof}

As a corollary, we obtain the following well-known result \cite{BGS}.
\begin{lemma} The metric $\hat{g}$ is extremal if and only if the induced metric $\check{g}$ is.
\end{lemma}

\begin{proof}
Fix a basis of $\hat{\mathfrak{t}}$, $e_0=b$, $e_1$, ... $e_n$, and identify $\hat{\mathfrak{t}}/(\bR b)$ with $\bR^n$ by using the basis $\{e_1, \dots, e_n\}$. Then, the $S^2\mathfrak{t}^*$--valued function of the quotient restricted to the hyperplane $\langle \hat{\mu},b\rangle=1$ becomes $\bfH^b = (H_{ij})=(\hat{H}_{ij} -\hat{\mu}_i\hat{\mu}_j)$. We then compute
$$s_{\hat{g}} = \sum_{i,j=0}^n \frac{\del^2}{\del \hat{\mu}_i \hat{\mu}_j}\hat{H}_{ij} = 0+ \sum_{i,j=1}^n \frac{\del^2}{\del \hat{\mu}_i \hat{\mu}_j}\hat{H}_{ij} =  \sum_{i,j=1}^n \frac{\del^2}{\del \hat{\mu}_i \hat{\mu}_j}H_{ij} + n^2 = s_{\check{g}} + n^2$$
by using Abreu's formula (\ref{abreuForm}).
\end{proof}

Combining Proposition~\ref{propHbSasak} with the result of~\cite{zz} we get the following corollary (where for simplicity we denote the labeled polytope $(\Delta_b, \varrho_b(\hat{u}_1),\dots , \varrho_b(\hat{u}_d))$ by $(\Delta_{b},u_{b})$).
\begin{corollary} Let $(N^{2n+1},g, \mD, \hat{T}^{n+1})$ be a compact toric Sasaki manifold with Reeb vector field $X_b$. If $g$ is extremal (in the sense that $\check{g}$ is) then $(\Delta_{b},u_{b})$ is analytically relatively $K$--stable with respect to toric degenerations. Moreover, if $g$ has constant scalar curvature then $\zeta_{(\Delta_{b},u_{b})}$ is constant.
\end{corollary}
We expect the converse to be true for $5$--dimensional contact toric manifolds with constant extremal affine function, in view of~\cite{don:scalar,don:estimate,don:extMcond,don:csc} but we do not show this in the present paper.

\section{Orthotoric structures and generic quadrilaterals.}\label{sectOrtho}
\subsection{Orthotoric structures}\label{sectGEN=ortho}\label{orthoDEFNandPARAM}

\begin{definition}\label{defnORTHO} Let $(M^4,\omega,J, g, T,\mu)$ be a compact, connected, K\"ahler toric\footnote{Every compact orbifold admitting an orthotoric K\"ahler metric in the sense of \cite{wsd} is toric, see~\cite{H2FI}.} $4$--orbifold. It is \defin{orthotoric} if there exist two positive smooth $T$--invariant functions $x$, $y\in C^{\infty}(\mathring{M})$ with $g$--orthogonal gradients on $\mathring{M}$ and an identification between $\mathfrak{t}^*$ and $\bR^2$ through which the moment map is $\mu=(x+y,xy)$. We call $x$, $y$ \defin{orthotoric coordinates} on $\mathring{M}$.
\end{definition}

In~\cite{H2FI}, it is shown that an orthotoric orbifold admits a Hamiltonian $2$--form of maximal order. Recall that a $2$--form $\Psi$ is Hamiltonian if, for any vector field $X$, $$2\nabla_X \Psi = d\sigma_1\wedge \omega(X,\cdot) + d^c\sigma_1 \wedge g(X,\cdot)$$ where
$\sigma_1 = tr_{\omega} \Psi = \frac{\Psi\wedge \omega}{\omega^2}.$ In complex dimension $2$, the orthotoric coordinates, $x$, $y$ are the eigenvalues of $J\circ\Psi$ (viewed as a field of complex endomorphisms of $(M,J)$ via the metric $g$). In particular, $x$ and $y$ are continuously defined on the whole $M$. Notice that $\sigma_1 = x+y$ and we set
$$\sigma_2 = xy = det_{\omega} \Psi = \frac{\Psi\wedge \Psi}{\omega^2}.$$ The moment map of an orthotoric orbifold is
$$\mu=(\sigma_1,\sigma_2)= (x+y,xy)$$
and its moment polytope $\Delta=\im \mu$ has a special shape which we shall now describe. Since $\mu$ is a moment map, it has rank $2$ on $\mathring{M}$ and so $x-y$ does not vanish on $\mathring{M}$. Our convention is that $x\geq y$ on $M$ and we set $\im x= [\alpha_1,\alpha_2]$ and $\im y = [\beta_1,\beta_2]$, with $\alpha_1\geq \beta_2$. The facets of $\Delta$ are explicitly given as the image via $\sigma \co(x,y) \mapsto (x+y,xy)$ of the facets of the rectangle $[\alpha_1,\alpha_2]\times[\beta_1,\beta_2]$, that is, for $i=1$, $2$:
$$F_{\alpha_i}=
\{ \sigma(\alpha_i,y) \,|\; y\in [\beta_1,\beta_2] \} \;\mbox{ and }\; F_{\beta_i}=\{ \sigma(x,\beta_i) \,|\; x\in [\alpha_1,\alpha_2] \}.$$ The $\Vlabel$s of $\Delta$ associated to $F_{\alpha_1}$, $F_{\alpha_2}$, $F_{\beta_1}$ and $F_{\beta_2}$ are, respectively
\begin{align} \label{labelOrtho} u_{\alpha_1}=C_{\alpha_1}\begin{pmatrix} \alpha_1 \\
-1
\end{pmatrix},\;\;u_{\alpha_2}=C_{\alpha_2}\begin{pmatrix} \alpha_2 \\
-1
\end{pmatrix}, \;\;u_{\beta_1}=C_{\beta_1}\begin{pmatrix} \beta_1 \\
-1
\end{pmatrix},\;\;u_{\beta_2}=C_{\beta_2}\begin{pmatrix} \beta_2 \\
-1
\end{pmatrix}
\end{align}
for some constants $C_{\alpha_1}$, $C_{\beta_2} >0$ and $C_{\alpha_2}$, $C_{\beta_1} <0$, where the signs are prescribed by the convention that the $\Vlabel$s are inward vectors. \begin{rem}\label{wProjspaceOrtho}
The case $\alpha_2=\beta_1$ has been extensively studied in~\cite{H2FII} and corresponds to triangles, that is weighted projective spaces. It is shown that the only $n$--dimensional compact manifold admitting an orthotoric structure is $\bcp^n$.
\end{rem}

\begin{definition}\label{defnOrthoPolyt}
A quadrilateral $\Delta \subset \bR^2$ is an \defin{orthotoric polytope} if there exists a rectangle $[\alpha_1,\alpha_2]\times[\beta_1,\beta_2]\subset \bR^2$, with $\beta_2<\alpha_1$, which is mapped on $\Delta$ by the map $\sigma(x,y)=(x+y,xy)$. We denote such a polytope $\Delta_{\alpha_1,\alpha_2,\beta_1,\beta_2}$. Thus, any labeled orthotoric polytope determines and is determined by $8$ real numbers $(\alpha_1,\alpha_2,\beta_1,\beta_2,C_{\alpha_1}, C_{\alpha_2}, C_{\beta_1},C_{\beta_2})$ which we shall refer to as \defin{orthotoric parameters}.
\end{definition}

Let $\Delta$ be a convex quadrilateral in $\mathfrak{t}^*$ and $s_1$, $s_2$, $s_3$, $s_4$ its vertices, such that $s_1$ and $s_3$ are not adjacent. By identifying $s_2-s_1$ and $s_4-s_1$ with a basis of $\mathfrak{t}^*$, $\Delta$ becomes the convex hull in $\bR^2$ of points $(0,0),(1,0),(a,b),(0,1)$. We denote the resulting polytope $\Delta'_{(a,b)}$ and call it a \defin{normal form} of $\Delta$. Since $\Delta$ is convex, $a>0$, $b>0$ and $a+b>1$. Notice that if $a$, $b \neq 1$, $\Delta$ is \defin{generic}, that is, has no parallel edges. Obviously, the normal form is not unique.

\begin{lemma} \label{lemNORMALrep}There exists an affine map of $\bR^2$ mapping $\Delta'_{(a',b')}$ to $\Delta'_{(a,b)}$ if and only if $(a',b')$ belongs to
\begin{align}\label{tabOrbS4}
\begin{split}
\bigg\{ & (a,b), \frac{1}{a}(a+b-1,1), \frac{1}{a+b-1}(a,b), \frac{1}{b}(1,a+b-1), \\
 & \;(b,a), \frac{1}{a}(1,a+b-1), \frac{1}{a+b-1}(a,b), \frac{1}{b}(a+b-1,1) \bigg\}.
\end{split}
\end{align}
In particular, a convex quadrilateral $\Delta$ is affinely equivalent to a unique $\Delta'_{(a,b)}\subset \bR^2$ such that $a$, $b>0$, $b\leq1 \leq a$ and $a + b \geq 2$.
\end{lemma}

\begin{lemma}\label{theoNoParall=ortho}
Any generic convex quadrilateral is equivalent to an orthotoric polytope.
\end{lemma}

\begin{proof}
Let $\Delta'_{(a,b)}$ be a normal form of a convex quadrilateral $\Delta$. The affine transform $$\Phi^* =\begin{pmatrix}
1-b &a-1\\
1-b&0\end{pmatrix} + \begin{pmatrix}
1\\
0\end{pmatrix}$$ has an inverse if and only if $\Delta$ is generic. In that case, $\Phi^*(\Delta'_{(a,b)})$ is the convex hull of $(1,0)$, $(a,0)$, $(2-b,(1-b))$ and $(a+(1-b),a(1-b))$. In particular, $\Phi^*(\Delta'_{(a,b)})$ is orthotoric, as the image $\sigma([\alpha_1,\alpha_2]\times[\beta_1,\beta_2])$, with
$$\alpha_1 =\min\{1,a\},\;\;\alpha_2 =\max\{1,a\},\;\;\beta_1 =\min\{0,1-b\},\;\;\beta_2 =\max\{0,1-b\}.$$
Notice that $\beta_2<\alpha_1$, since $a+b>1$ and $a$, $b>0$.
\end{proof}

\begin{corollary} \label{coroCharactPair} For any convex quadrilateral $\Delta$, there exists a unique pair $(\alpha,\beta)$, satisfying $0\leq\beta<1\leq\alpha$ and $\alpha-\beta \geq 1$, such that $\Delta$ is affinely equivalent to the normal form $\Delta'_{(\alpha,1-\beta)}$. We call $(\alpha,\beta)$ the \defin{characteristic pair} of $\Delta$. If $\Delta$ is generic then $\beta>0$ (and thus $\alpha>1$) and it is affinely equivalent to the orthotoric quadrilateral $\Delta_{\alpha,\beta}= \sigma([0,\beta]\times[1,\alpha])$.
\end{corollary}

We now summarize and rephrase~\cite[Proposition 1]{H2FII} and~\cite[Proposition 11]{H2FI} which together describe the space of orthotoric K\"ahler metrics compatible with a given toric symplectic orbifold.
\begin{proposition}\label{PropOrthometriclocale}\cite{H2FII},\cite{H2FI} Let $(M,\omega,g,\mu,T)$ be an orthotoric orbifold with orthotoric coordinates $x$ and $y$ and momentum coordinates $\sigma_1=x+y$, $\sigma_2=xy$. Let $t_1$, $t_2$ be the corresponding angle coordinates on $\mathring{M}$. Letting $\im x =[\alpha_1,\alpha_2]$ and $\im y = [\beta_1,\beta_2]$, there exist functions, $A\in C^{\infty}([\alpha_1,\alpha_2])$ and $B\in C^{\infty}([\beta_1,\beta_2])$, such that $A(x)$ and $B(y)$ are positive on $\mathring{M}$,  \begin{equation}\label{metric}
g_{|_{\mathring{M}}}= \frac{(x-y)}{A(x)}dx^2 +\frac{(x-y)}{B(y)}dy^2 + \frac{A(x)}{x-y}(dt_1+ y dt_2)^2 + \frac{B(y)}{x-y}(dt_1+ x dt_2)^2
\end{equation} and \begin{equation}\label{C1}\begin{split}
   A(\alpha_i)=0, \;\; &B(\beta_i)=0,\\
   A'(\alpha_i) = 2/C_{\alpha_i},\;\; &B'(\beta_i) = -2/C_{\beta_i}.
\end{split} \end{equation}

Conversely, for any smooth functions respectively positive on $(\alpha_1,\alpha_2)$ and $(\beta_1,\beta_2)$ and satisfying~(\ref{C1}), the formula~(\ref{metric}) defines a smooth orthotoric K\"ahler metric on $M$ compatible with $\omega$, with orthotoric coordinates $x$ and $y$.
\end{proposition}

The $S^2\bR^2$--valued function $\bfH=(g(X_{e_i}, X_{e_j}))$ associated to the metric~(\ref{metric}) is
\begin{equation}\label{bfHABortho}
\bfH_{A,B} =\frac{1}{x-y}
\begin{pmatrix} A(x) + B(y) & yA(x)+xB(y)\\
 yA(x)+xB(y)& y^2A(x)+x^2B(y)
\end{pmatrix}.
\end{equation}
Proposition~\ref{PropOrthometriclocale} can be derived from~(\ref{bfHABortho}) by using Proposition~\ref{prop1H2FII}. The integrability condition (\ref{integCond}) is satisfied by construction. In particular, one can give the explicit symplectic potential~\cite[Proposition 11]{H2FI}:
\begin{equation} \label{POTortho}
G_{A,B}(x,y) =-\int_{\alpha_1}^x \frac{(t-x)(t-y)}{A(t)} dt + \int_{\beta_1}^y \frac{(t-x)(t-y)}{B(t)}dt.
\end{equation}
From Lemma~\ref{theoNoParall=ortho} and Proposition~\ref{PropOrthometriclocale}, we infer.
\begin{proposition}\label{coroOrtPolyt=orthoMet}\label{theoGENquad=ortho}
Let $(M,\omega,T,\mu)$ be a symplectic toric $4$--orbifold. There exists an orthotoric metric compatible with $\omega$ if and only if the polytope $\im \mu$ is a generic quadrilateral or a triangle.
\end{proposition}
\begin{proof} For triangles, see~\cite{H2FII}. An orthotoric quadrilateral is generic by definition. Conversely, if $\im \mu$ is generic, by Lemma~\ref{theoNoParall=ortho}, we can identify the Lie algebra of $T$ with $\bR^2$, such that $\im \mu$ be an orthotoric polytope. Thus, $\im \mu$ is the image by $\sigma$ of a rectangle $[\alpha_1,\alpha_2]\times[\beta_1,\beta_2]\subset \bR^2$ with $\beta_2<\alpha_1$. The functions $x= \frac{\sigma_1+\sqrt{\sigma_1^2- 4\sigma_2}}{2\sigma_2}$ and $y= \frac{\sigma_1-\sqrt{\sigma_1^2- 4\sigma_2}}{2\sigma_2}$ are smooth on $M$ since $x>y$. Moreover, $dx$ and $dy$ are linearly independent on $\mathring{M}$ since $\mu$ is a moment map. Hence, taking positive functions satisfying condition~(\ref{C1}) leads to a smooth metric~(\ref{metric}) for which $dx$ and $dy$ are orthogonal.
\end{proof}

\subsection{Extremal orthotoric metrics}\label{sectOrthometric}\label{sectCPTextAB}\label{sectCSCsolve}\label{sectEXTsolve}

The ``separation of variables'', mentioned in the introduction, appears in the formula giving the scalar curvature of an orthotoric metric $g_{A,B}$, \cite{wsd}:
\begin{equation}\label{curvature}
S(G_{g_{A,B}}) = -\frac{A^{''}(x) + B^{''}(y)}{x-y}.
\end{equation}
\begin{rem}\label{remReferee} It is elementary to verify that if $A$ and $B$ are positive, the inverse Hessian of the potential~(\ref{POTortho}), with respect to $\sigma_1=x+y$ and $\sigma_2=xy$, is~(\ref{bfHABortho}) and then that the scalar curvature is~(\ref{curvature}). It is not necessary to notice the presence of a Hamiltonian $2$--form (orthotoric metric in this case) but it gives a unified and geometric framework.
\end{rem}
 Then, the condition for the metric to be extremal (which, in the toric context, amounts to the fact that its scalar curvature is the pull back by the moment map of an affine-linear function on $\mathfrak{t}^*$ with respect to variables $\sigma_1=x+y$ and $\sigma_2=xy$) gives rise to the following conditions on $A$ and $B$:
\begin{proposition}~\cite{wsd}\label{condEXTREMprop}
 Let $(M,\omega,g,J,\mu)$ be an orthotoric orbifold where $g$ is expressed on $\mathring{M}$ as~(\ref{metric}) with respect to orthotoric coordinates $x$,$y$. Then, $g$ is extremal if and only if $A$ and $B$ are polynomials of degree $4$, say $A(x) = A_0x^4 + A_1x^3 + A_2x^2 + A_3x + A_4$ and $B(y) = B_0y^4 + B_1y^3 + B_2y^2 + B_3y + B_4$, with
\begin{equation}
\label{condEXTREM}
 A_0 = -B_0, \;\; A_1 = -B_1 \;\mbox{ and }\; A_2 = -B_2.
\end{equation}Moreover, assuming~(\ref{condEXTREM}), $g$ has constant scalar curvature if and only if $A_0=0$ and is K\"ahler--Einstein if and only if
\begin{equation}\label{KE}
A_0=0 \mbox{ and } A_3=-B_3.
\end{equation}
\end{proposition}

\begin{rem}\label{sigma1parall} Notice that for any extremal orthotoric metric \begin{equation}
  S(G_{g_{A,B}}) = -12A_0(x+y) -6A_1 = -12A_0\sigma_1 -6A_1.
\end{equation}
It follows that the scalar curvature of an extremal orthotoric metric is equipoised in the sense of Definition~\ref{invOrthoFct}.
\end{rem}
\begin{rem}\label{abstractProb} Although we found it geometrically instructive to present most of the material of this section in the setting of toric orbifolds, the assumption that the corresponding labeled polytopes are rational is unnecessary for the results to hold true. For instance, if $(\Delta,u)$ is a labeled orthotoric polytope with orthotoric parameters $(\beta_1,\beta_2,\alpha_1,\alpha_2,C_{\alpha_1}, C_{\alpha_2}, C_{\beta_1}, C_{\beta_1})$ and if there exist polynomials $A$ and $B$ respectively positive on $[\alpha_1,\alpha_2]$ and $[\beta_1,\beta_2]$, satisfying conditions~(\ref{C1}) and~(\ref{condEXTREM}), then the function $G$ defined by (\ref{POTortho}) is a solution of Problem~\ref{problem1} for $(\Delta,u)$. Moreover, via Propositions~\ref{PropOrthometriclocale} and~\ref{condEXTREMprop}, the metric $g_{A,B}$ given by~(\ref{metric}) defines an extremal K\"ahler metric on $\mathring{\Delta}\times \bR^n$ with boundary condition given by~(\ref{condCOMPACTIFonH}).
\end{rem}

\noindent The next lemma states a condition on the orthotoric parameters for that it exists a \defin{formal solution} of Problem~\ref{problem1}, that is a solution $\bfH_{A,B}$ of equation~(\ref{extrem}) which is not necessarily positive definite. This is equivalent to the fact that there exist polynomials $A$ and $B$ satisfying conditions~(\ref{C1}) and~(\ref{condEXTREM}) but they are NOT necessarily respectively positive on $[\alpha_1,\alpha_2]$ and $[\beta_1,\beta_2]$.

\begin{lemma}\label{condEXTlemma} Let $\beta_1<\beta_2<\alpha_1<\alpha_2$ be real numbers and $C_{\alpha_1}$, $C_{\alpha_2}$, $C_{\beta_1}$, $C_{\beta_1}$ be non-zero real numbers. There exist polynomials of degree at most $4$, $A$ and $B$, satisfying conditions~(\ref{C1}) and~(\ref{condEXTREM}) if and only if
\begin{equation}
\label{condEXTpoly5}
\begin{split}
&\frac{(\alpha_1 +\alpha_2)^2+2\alpha^2_2 +(\beta_1 +\beta_2)^2 +2\beta_1\beta_2 -2(2\alpha_2 + \alpha_1)(\beta_1 +\beta_2)}{(\alpha_2-\alpha_1)^2 \, C_{\alpha_1}} \\
+& \frac{(\alpha_1 +\alpha_2)^2+2\alpha^2_1 +(\beta_1 +\beta_2)^2 +2\beta_1\beta_2 -2(2\alpha_1 + \alpha_2)(\beta_1 +\beta_2) }{(\alpha_2-\alpha_1)^2 \, C_{\alpha_2}}\\
+&\frac{(\beta_1 +\beta_2)^2+2\beta^2_2 +(\alpha_1 +\alpha_2)^2 +2\alpha_1\alpha_2 -2(2\beta_2 + \beta_1)(\alpha_1 +\alpha_2) }{(\beta_2-\beta_1)^2 \, C_{\beta_1} }\\
+&\frac{(\beta_1 +\beta_2)^2+2\beta^2_1 +(\alpha_1 +\alpha_2)^2 +2\alpha_1\alpha_2 -2(2\beta_1 + \beta_2)(\alpha_1 +\alpha_2)}{(\beta_2-\beta_1)^2 \, C_{\beta_2}} =0.\\
\end{split}
\end{equation} In this case, the polynomials $A$ and $B$ are uniquely determined by the orthotoric parameters. They are of degree $3$ (i.e $A_0=0$) if and only if \begin{align}
\frac{1}{(\alpha_2-\alpha_1)^2}\left(\frac{1}{C_{\alpha_1}} + \frac{1}{C_{\alpha_2}} \right) &= \frac{1}{(\beta_2-\beta_1)^2} \left(\frac{1}{C_{\beta_1}} + \frac{1}{C_{\beta_2}}\right).\label{anciennement-ii}
\end{align} $A$ and $B$ satisfy the condition (\ref{KE}) if, in addition to~(\ref{condEXTpoly5}) and~(\ref{anciennement-ii}),
\begin{align}\label{anciennement-iv}
\begin{split}
\frac{1}{(\alpha_2-\alpha_1)^2} & \left(\frac{\alpha_2(2\alpha_1 + \alpha_2)}{C_{\alpha_1}} + \frac{\alpha_1(2\alpha_2 + \alpha_1)}{C_{\alpha_2}} \right)\\
& = \frac{1}{(\beta_2-\beta_1)^2} \left(\frac{\beta_2(2\beta_1 + \beta_2)}{C_{\beta_1}} + \frac{\beta_1(2\beta_2 + \beta_1)}{C_{\beta_2}}\right).
\end{split}
\end{align}
\end{lemma}

\begin{proof} The compactification condition~(\ref{C1}) implies that $\alpha_1$, $\alpha_2$ and $\beta_1$, $\beta_2$ are roots of $A$ and $B$ respectively. So if such polynomials exist, they must be of the form
\begin{equation}
\label{polyEXT}
\begin{split}
A(x) &= (x-\alpha_1)(x-\alpha_2)(A_0x^2+R_1x +R_2), \\
B(y) &= (y-\beta_1)(y-\beta_2)(-A_0y^2+S_1y +S_2),
\end{split}
\end{equation} for some constant $A_0, R_1, R_2, S_1, S_2$. The extremality conditions~(\ref{condEXTREM}) imply
\begin{equation}
\label{A1=-B1}
\begin{split}
R_1 -A_0(\alpha_1 +\alpha_2) &= -S_1 -A_0 (\beta_1 +\beta_2), \\
R_2 +A_0\alpha_1\alpha_2 - R_1(\alpha_1 +\alpha_2) &= -S_2 + A_0 \beta_1\beta_2 + S_1(\beta_1 +\beta_2).
\end{split}
\end{equation}
 The compactification condition~(\ref{C1}) reads, in terms of~(\ref{polyEXT}),
\begin{equation}
\label{condDeriveeAi}
\begin{split}
A'(\alpha_i) = (-1)^i(\alpha_2-\alpha_1) (A_0\alpha_i^2+R_1\alpha_i +R_2)= 2/C_{\alpha_i},\\
B'(\beta_i) =(-1)^i(\beta_2-\beta_1) (-A_0\beta_i^2+S_1\beta_i +S_2)=-2/C_{\beta_i},
\end{split}
\end{equation}for $i=1$,$2$. $R_1$ and $R_2$ may be expressed as functions of $A_0$, by using the condition~(\ref{condDeriveeAi})
\begin{equation}\label{definXI}
\begin{split}
R_1 &= \frac{1}{(\alpha_2-\alpha_1)^2}\left(\frac{2}{C_{\alpha_1}} + \frac{2}{C_{\alpha_2}} \right) - A_0 (\alpha_1 +\alpha_2),\\
R_2 &= -\frac{1}{(\alpha_2-\alpha_1)^2}\left(\frac{2\alpha_2}{C_{\alpha_1}} + \frac{2\alpha_1}{C_{\alpha_2}} \right) + A_0\alpha_1\alpha_2,
\end{split}
\end{equation}
and similarly,
\begin{equation}\label{definETA}
\begin{split}
S_1 &= -\frac{1}{(\beta_2-\beta_1)^2} \left(\frac{2}{C_{\beta_1}} + \frac{2}{C_{\beta_2}}\right)+A_0(\beta_1 +\beta_2),\\
S_2 &= \frac{1}{(\beta_2-\beta_1)^2} \left(\frac{2\beta_2}{C_{\beta_1}} + \frac{2\beta_1}{C_{\beta_2}}\right)- A_0\beta_1\beta_2.
\end{split}
\end{equation} The condition~(\ref{A1=-B1}), together with the fact that $\alpha_1 +\alpha_2- (\beta_1 +\beta_2)\neq 0$, allows us to determine
\begin{equation}\label{definA0}
\begin{split}
 A_0 = \frac{1}{2(\alpha_1 +\alpha_2- (\beta_1 +\beta_2))}\bigg(\frac{1}{(\alpha_2-\alpha_1)^2}\left(\frac{2}{C_{\alpha_1}} + \frac{2}{C_{\alpha_2}} \right) \qquad \qquad \\
-\frac{1}{(\beta_2-\beta_1)^2} \left(\frac{2}{C_{\beta_1}} + \frac{2}{C_{\beta_2}}\right) \bigg).
\end{split}
\end{equation} The second line of~(\ref{A1=-B1}) implies
\begin{equation}
\label{condEXTpoly}
\begin{split}
\frac{2}{(\alpha_2-\alpha_1)^2}\left(\frac{2\alpha_2 + \alpha_1}{C_{\alpha_1}} + \frac{2\alpha_1 + \alpha_2}{C_{\alpha_2}}\right)  - \frac{2}{(\beta_2-\beta_1)^2} \left(\frac{2\beta_2+\beta_1}{C_{\beta_1}} + \frac{2\beta_1+\beta_2}{C_{\beta_2}}\right)  \\
= A_0 \big((\alpha_1+ \alpha_2)^2 + 2\alpha_1\alpha_2 -(\beta_1+\beta_2)^2-2\beta_1\beta_2\big)
\end{split}
\end{equation} which may be written as~(\ref{condEXTpoly5}). Equation~(\ref{anciennement-ii}) comes directly from the line~(\ref{definA0}). Finally, the condition~(\ref{KE}) is $A_3=-B_3$ assuming $A_0=0$, that is $-R_1 \alpha_1\alpha_2+ R_2(\alpha_1+\alpha_2) = -S_1 \beta_1\beta_2+ S_2(\beta_1+\beta_2)$ is exactly equation~(\ref{anciennement-ii}).
\end{proof}
\begin{rem}\label{positivityC+} Polynomials of degree $3$ satisfying the compactification condition~(\ref{C1}) are automatically positive on $(\alpha_1,\alpha_2)$, and $(\beta_1,\beta_2)$.
\end{rem}

\section{Calabi toric structures and trapezoids.}\label{sectCalabi}
\subsection{K\"ahler toric $4$-orbifolds with non-trivial Hamiltonian $2$--forms.}

\begin{definition} \label{defnCALABItoric1}A K\"ahler toric $4$--orbifold, $(M,\omega,J,g, T,\mu)$, is \defin{Calabi toric} with respect to a Killing vector field $K$, if there exists a Hamiltonian $2$--form of order $1$ whose non-constant eigenvalue is a Hamiltonian function of $K$.\end{definition}

From \cite[Proposition 10]{H2FII} we know that, excluding $\bcp^2$, there is a unique non trivial Hamiltonian $2$--form up to addition of a multiple of the symplectic form, on a toric $4$--orbifold. In particular, such Hamiltonian $2$--form, as well as the symplectic gradient of its trace $K$, is $T$--invariant.
\begin{proposition} \label{propHam2fOrder1polytope}
 The moment polytope of a symplectic toric orbifold admitting a compatible Calabi toric metric is a trapezoid or a triangle.\end{proposition}
\begin{proof}Suppose that $(M,\omega,J, g, K, T,\mu)$ is Calabi toric. Let $v\in \mathfrak{t}$ be such that $K=X_{v}$, we know that such vector $v$ exists since $K$ is a Hamiltonian Killing vector field commuting with the infinitesimal action of $T$ on $M$. Let $x=\mu(v)$. Notice that $dx =d\mu(v)=-\omega(X_v,\cdot)$ vanishes only on zero set of $K$ while, thanks to \cite[Theorem 1]{H2FI}, we know that $g(K,K)$ only depends on the value $x$. In particular, critical points of $x$ are only contained in the preimage of the ends of the interval $\im \,x$. Hence, $\Delta\cap \mbox{ev}_v^{-1}(\mathring{\im\,x})$ does not contain any vertex of $\Delta$. Thus, $\Delta$ is a trapezoid or a triangle.\end{proof}

From now on, we assume that the toric action of a Calabi toric orbifold fixes four points, so its polytope is a trapezoid and we exclude weighted projective space of our consideration.

\begin{proposition}\label{lemCALABItoricFunc} The K\"ahler toric orbifold $(M,\omega,J,g, T,\mu)$ is Calabi toric with respect to $K$ if and only if there exist two positive smooth $T$--invariant functions $x,y\in C^{\infty}(M)$ with $g$--orthogonal gradients such that $K= J\mbox{grad }x$ and there is an identification between $\mathfrak{t}^*$ and $\bR^2$ for which $$\mu=(x,xy).$$\end{proposition}

\begin{proof} Suppose first that $(M,\omega,J, g, K, T,\mu)$ is Calabi toric and consider the same notation then the one in the proof of Proposition~\ref{propHam2fOrder1polytope}. Recall that $\mathring{M}$ denotes the open dense subset where the torus acts freely. By the general theory, $\mathring{M}$ is a manifold diffeomorphic to $\mathring{\Delta}\times T$. Hence, the manifold $\mathring{M}$, when endowed with the restriction of the K\"ahler structure of $M$, is a connected K\"ahler manifold admitting a Hamiltonian $2$--form of order one whose only non-constant eigenvalue, say $\xi$, is a Killing potential for $K$. This situation is exactly the one required to apply \cite[Theorem 1]{H2FI}. Thus, we obtain the explicit description of the restriction of the K\"ahler structure $(\omega,J, g)$ on $\mathring{M}$:
\begin{equation}\label{metCalabiWSD}\begin{split}
   g = \frac{x}{A(x)}dx^2 + x g_{\lambda} + \frac{A(x)}{x} &\theta^2, \;\;\;\;\;Jdx = \frac{A(x)}{x} \theta, \\
   \omega= x\omega_{\lambda} + dx\wedge\theta, \;\;\;\;\; &\;\;\;\;\; d\theta =-\omega_{\lambda}
 \end{split}\end{equation}
where $\lambda$ is the constant root of the Hamiltonian $2$-form so that $x:=|\lambda-\xi|$ is positive on $\mathring{M}$, the $1$--form $\theta$ satisfies $\theta(K)=1$ and $(g_{\lambda},\omega_{\lambda})$ is a K\"ahler structure on a Rieman surface.

Adding a constant to $\mu$ if necessary, we suppose that $x=\mu(v)>0$. Let $u\in \mathfrak{t}$ be any vector for which the pair $(v,u)$ is a basis of $\mathfrak{t}$. This choice gives an identification between $\mbox{Lie } T$ and $\bR^2$, so that $$\mu=\mu(v) v^* + \mu(u) u^*= (x,\mu(u)).$$

Translating by $v$ if necessary, one can choose $u$ in order that the $T$-invariant function $$f_u=\theta(X_u)= \frac{x}{A(x)}g(X_u,X_v)$$ be non-negative. One has $\mL_{X_u} f_u =0$, $\mL_{X_v} f_u =0$ from $T$--invariance and \begin{equation}\label{differf_ucalabi}\begin{split}
  \mL_{JX_v}f_u &= d\theta (JX_v,X_u) =0\\
  \mL_{JX_u}f_u &= -d\theta(X_u,JX_u) =\frac{-1}{x}(\omega(X_u,JX_u) + f_udx(JX_u))
\end{split}
\end{equation} using~(\ref{metCalabiWSD}). The vector fields $X_u,X_v, JX_u,JX_v$ commute so the formulae above determine the differential $df_u =\frac{1}{x} (d\mu(u)- f_udx)$ and then $\mu(u)= xf_u+ c$ for a certain constant $c$.

The map $\tilde{\mu}=\mu-cu^*$ is a moment map for the action of $T$ such that $\tilde{\mu}(v)=\mu(v)>0$ and \begin{equation}\label{momentCalabi}\tilde{\mu}= \mu(v) v^* + (\mu(u)-c) u^*= (x,xf_u).\end{equation} In particular, the function $f_u$ may be defined on the whole $M$. Indeed, the assumption that $T$ fixes $4$ points of $M$ is equivalent to that $\im \tilde{\mu}$ is a quadrilateral which can be true only if $x>0$ on $M$. We set $y=f_u$, so $\tilde{\mu}=(x,xy)$. Finally, notice that the first formula of~(\ref{differf_ucalabi}) implies that $dy$ and $dx$ are orthogonal.

Conversely, the fact that the moment map is $\mu=(x,xy)$, with $x>0$ and $y\geq0$, implies that the rational labeled polytope associated to $(M,\omega,T)$ has a $\Vlabel$ lying in $\bR e_1$. In particular, the line $\bR e_1$ meets $\Lambda$ (where $T=\bR^2/\Lambda$) and determines a circle, $S\subset T$, the generator of which is $\bR$-collinear to $K=J\mbox{grad} x$ and is a Hamiltonian Killing vector field on $M$. Consider the set of action-angle coordinates $(\sigma_1,\sigma_2,t,s)$ on $\mathring{M}$ where $(\sigma_1,\sigma_2)=(x,xy)$ and $dt(K)=1$ and the normal expression for the toric metric in term of matrices $\bfH=(H_{ij})$ and $\bfG=(G_{ij})$. The assumption that $dx$ and $dy$ are $g$--orthogonal implies the following relations:
\begin{equation}\label{calabiOrthogRelat}
  G_{11} + 2y G_{12} + y^2G_{22} =f(x,y), \;\; x(G_{12} + yG_{22}) = 0 \;\mbox{ and }\; x^2G_{22} =h(x,y).
\end{equation}

From these equations and the integrability of $J$, we infer \begin{equation}\begin{split}
 &\frac{\del}{\del \sigma_2} G_{11} - \frac{\del}{\del \sigma_1} G_{12} =\frac{1}{x}\del_{y}f - y(\frac{\del}{\del \sigma_2} G_{12} - \frac{\del}{\del \sigma_1} G_{22}) \equiv0, \\
 &\frac{\del}{\del \sigma_2} G_{12} - \frac{\del}{\del \sigma_1} G_{22} =\frac{1}{x^3}(h-x\del_x h)\equiv0.
 \end{split}
 \end{equation}
Hence, there exist functions of one variable $A$ and $B$ such that $f(x,y)=x/A(x)$, $h(x,y) = x/B(y)$ and
\begin{align} \bfH_{A,B} = \frac{1}{x}\begin{pmatrix} A(x) & yA(x)\\
yA(x) & x^2B(y) + y^2A(x)
\end{pmatrix},\;\; \bfG_{A,B} = \begin{pmatrix} \frac{x}{A(x)} +\frac{y^2}{xB(y)}& \frac{-y}{xB(y)} \\ \frac{-y}{xB(y)} & \frac{1}{xB(y)}
\end{pmatrix}.\label{bfHABcalabi}
\end{align} Hence, we obtain that the expression of the K\"ahler structure on $\mathring{M}$ associated to $\bfH_{A,B}$ and $\bfG_{A,B}$ is exactly of the form~(\ref{metCalabiWSD}) where $\theta=dt+yds$, $\omega_{\lambda} = dy\wedge ds$ and $g_{\lambda}=\frac{dy^2}{B(y)} +B(y)ds^2$. Using~\cite[Theorem 2]{H2FI}, one gets a Hamiltonian $2$--form of order one on $\mathring{M}$ given explicitly by $\Psi = x(\omega-x\omega_{\lambda})$. This $2$--form admits a unique extension on $M$ since it is a parallel $2$--form with respect to a connection~\cite[Proposition 4]{H2FI}. \end{proof}

The classification of Hamiltonian $2$--forms on symplectic toric $4$--orbifolds claimed in the introduction follows from the following proposition together with Proposition~\ref{lemCALABItoricFunc} and the fact that orthotoric metrics characterize the ones admitting Hamiltonian $2$-forms of order two~\cite{H2FI}.
\begin{proposition} \label{propH2Forder0} Let $(M,\omega,J,g, T)$ be a K\"ahler toric $4$--orbifold admitting a non-trivial Hamiltonian $2$--form of order zero. Then, the moment polytope is a parallelogram. In particular, $(M,\omega,J,g)$ is a product of weighted projective lines endowed with a product of toric K\"ahler metrics.
\end{proposition}
\begin{proof} Let $\Psi$ be a non-trivial Hamiltonian $2$--form of order zero on $M$. Seen as an endomorphism using $g$, $J\Psi$ has two constant roots, say $\lambda_1$, $\lambda_2$. These roots are distinct since $\Psi$ is non-trivial (that is, it is not a multiple of the identity). The K\"ahler structure restricts to a K\"ahler structure on the eigenspaces of $\Psi$ and the induced splitting of the tangent space $TM=V_1\oplus V_2$ is invariant by the local action generated by $\mathfrak{t}$. Since $\Psi$ is parallel with respect to the Levi-Civita connection, the distributions $V_i$ are not only integrable in the Frobenius sense but also closed for the Levi-Civita connection meaning that $\nabla_Y X(p)\in (V_i)_p$ as soon as $X$ is a section of $V_i$. Thus, a vector field, $K$, such that $K=K_1 +K_2$ with $K_i\in V_i$, is Killing if and only if each vector field $K_i$ is. Hence, one can choose a basis $(e_1,e_2)$ of $\mathfrak{t}$ such that $g(X_{e_1},X_{e_2})=0$ and $X_{e_i}(p)\in (V_i)_p$ for all $p\in U$. Taking any moment map, $\mu$, of the action of $T$ on $M$, the functions $x=\mu(e_1)$ and $y=\mu(e_2)$ have orthogonal gradients which are non-vanishing on $\mathring{M}$. With respect to the basis $(e_1,e_2)$, $\im \mu$ is a rectangle since $\mu = (x,y).$ The Proposition follows from the Delzant--Lerman--Tolman Classification of toric orbifolds.
\end{proof}

\subsection{Calabi toric structures.}\label{sectCalabiPolyt}

In view of Proposition~\ref{lemCALABItoricFunc}, an equivalent definition of Calabi toric orbifold is
\begin{definition} \label{defnCALABItoric}
Let $(M,\omega,J,g, T,\mu)$ be a compact, connected, K\"ahler toric $4$--orbifold. It is \defin{Calabi toric} if there exist smooth $T$--invariant functions $x$ and $y\in C^{\infty}(M)$ with $x>0$, $y>0$ $g$--orthogonal gradients on $\mathring{M}$ and an identification between $\mathfrak{t}^*$ and $\bR^2$ through which the moment map is $\mu=(x,xy)$. We call $x$,$y$ the \defin{Calabi coordinates}.
\end{definition}

The moment polytope $\Delta= \im \mu$ of a Calabi toric orbifold has a special shape. Let $\im x = [\alpha_1,\alpha_2]$ and $\im y =[\beta_1,\beta_2]$, with $\alpha_1>0$ and $\beta_1\geq 0$. $\Delta$ is the image of $[\alpha_1,\alpha_2]\times[\beta_1,\beta_2]$, by $\sigma \co (x,y) \mapsto (x,xy)$. The $\Vlabel$s of $\Delta$ are
\begin{align} \label{labelCalabi}u_{\alpha_1}=C_{\alpha_1}\begin{pmatrix}
\alpha_1\\
0\end{pmatrix}, \; u_{\alpha_2}=C_{\alpha_2}\begin{pmatrix}
\alpha_2\\
0\end{pmatrix}, \;u_{\beta_1}=C_{\beta_1}\begin{pmatrix}
\beta_1\\
-1\end{pmatrix},\; u_{\beta_2}=C_{\beta_2}\begin{pmatrix}
\beta_2\\
-1\end{pmatrix}\end{align}
with $C_{\alpha_1}$, $C_{\beta_2}>0$ and $C_{\alpha_2}$, $C_{\beta_1}<0$. This leads us to the following definition.

\begin{definition} \label{defnCalabiPolyt}
A \defin{Calabi polytope} is a polytope which is the image of a rectangle $[\alpha_1,\alpha_2]\times[\beta_1,\beta_2]\subset\bR^2$, with $\alpha_1>0$ and $\beta_1\geq0$, by the map $\sigma \co (x,y) \mapsto (x,xy)$. Thus, any labeled Calabi polytope determines and is determined by a $8$--tuple $(\alpha_1,\alpha_2,\beta_1,\beta_2,C_{\alpha_1}, C_{\alpha_2}, C_{\beta_1},C_{\beta_2})$ we shall refer to as \defin{Calabi parameters}.
\end{definition}

\begin{lemma}\label{lemTrap=Calab}
A Calabi polytope is a trapezoid and any trapezoid is equivalent to a Calabi polytope.
\end{lemma}
\begin{proof}
Let $\alpha>1$ and $\Delta$ be the Calabi trapezoid given as the image by $\sigma$ of $[1,\alpha]\times[0,1]$, the affine map
$$\begin{pmatrix}
  1/(\alpha-1) & 0\\
  0& 1
\end{pmatrix} -\begin{pmatrix}
  1/(\alpha-1) \\
  0
\end{pmatrix}$$
maps $\Delta$ on $\Delta'_{(\alpha,1)}$, the convex hull in $\bR^2$ of $(0,0)$, $(1,0)$, $(1,\alpha)$, $(0,1)$. We conclude the proof, by using normal forms of quadrilaterals of Corollary~\ref{coroCharactPair}.
\end{proof}
\begin{rem} \label{CalabiNormForm}The Calabi trapezoid of parameters
$$(\alpha_1,\alpha_2,\beta_1,\beta_2,C_{\alpha_1}, C_{\alpha_2}, C_{\beta_1},C_{\beta_2})$$
is equivalent to the Calabi trapezoid of parameters $$\left(1,\alpha=\frac{\alpha_2}{\alpha_1},0,1,\alpha_1^2C_{\alpha_1}, \alpha_1^2C_{\alpha_2},(\beta_2-\beta_1) \alpha_1C_{\beta_1},(\beta_2-\beta_1) \alpha_1C_{\beta_1}\right).$$
\end{rem}
\begin{proposition} \label{defnCALABIlocalMet} Let $(M,\omega,g,\mu,T)$ be a Calabi toric orbifold with Calabi coordinates $x$,$y$ and momentum coordinates $\sigma_1=x$, $\sigma_2=xy$. Let $t_1$, $t_2$ be the corresponding angle coordinates on $\mathring{M}$. Letting $\im x =[\alpha_1,\alpha_2]$ and $\im y = [\beta_1,\beta_2]$, there exist functions, $A\in C^{\infty}([\alpha_1,\alpha_2])$ and $B\in C^{\infty}([\beta_1,\beta_2])$, such that $A(x)$ and $B(y)$ are positive on $\mathring{M}$,\begin{align}
g_{|_{\mathring{M}}} = x\frac{dx^2}{A(x)} + & x\frac{dy^2}{B(y)}  + \frac{A(x)}{x}(dt + yds)^2 + xB(y)ds^2 \label{CALABItoricmetric}
\end{align} and \begin{align}\label{eq:CondCompactCalabi}
\begin{split}
A(\alpha_i)=0, &\;\; B(\beta_i)=0 \\
A'(\alpha_i) =\frac{2}{C_{\alpha_i}}, &\;\; B'(\beta_i) =-\frac{2}{C_{\beta_i}}.
\end{split}
\end{align}

Conversely, for any smooth functions, $A$, $B$, respectively positive on $(\alpha_1,\alpha_2)$ and $(\beta_1,\beta_2)$ and satisfying~(\ref{eq:CondCompactCalabi}), the formula~(\ref{CALABItoricmetric}) defines a smooth Calabi toric metric on $M$ compatible with $\omega$, with Calabi coordinates $x$, $y$.
\end{proposition}

\begin{proof} The main part of the proof is a corollary of Proposition~\ref{lemCALABItoricFunc}. We just have to verify that the compactification conditions (\ref{condCOMPACTIFonH}) on $\bfH_{A,B}$ correspond to~(\ref{eq:CondCompactCalabi}): Then the proposition would follow from Proposition~\ref{prop1H2FII}. We compute that
$$\bfH_{(\alpha_i,y\alpha_i)}(u_{\alpha_i},\cdot) = C_{\alpha_i} \begin{pmatrix} A(\alpha_i)\\
yA(\alpha_i)\end{pmatrix} \;\;\mbox{and}\;\; \bfH_{(x,x\beta_i)}(u_{\beta_i},\cdot) = C_{\beta_i} \begin{pmatrix} 0\\
xB(\beta_i)\end{pmatrix}.$$
Moreover, $d\bfH_{(\alpha_i,y\alpha_i)}(u_{\alpha_i},u_{\alpha_i}) =  C_{\alpha_i}^2 A'(\alpha_i) \alpha_id\sigma_1$ and $$d\bfH_{(x,x\beta_i)}(u_{\beta_i},u_{\beta_i}) = -C_{\beta_i}^2 B'(\beta_i)(\beta_id\sigma_1 - d\sigma_2).$$

Finally, we have $u_{\beta_i} =C_{\beta_i} (\beta_id\sigma_1 - d\sigma_2)$ and $u_{\alpha_i}=C_{\alpha_i} \alpha_i d\sigma_1$ via the standard identifications $T_{\mu}\Delta \simeq \mathfrak{t}^*$ and $\mathfrak{t}^*\simeq \mathfrak{t}$, given by the basis $(\sigma_1,\sigma_2)$.
\end{proof}

\begin{corollary}
A symplectic potential of the Calabi toric metric $g_{A,B}$ is given by $$G(x,y) =x \int^{y}_{s_0} \left(\int^{s}_{t_0} \frac{1}{B(t)} dt \right)ds - x \int^{x}_{s_0} \frac{1}{s^2}\left(\int^{s}_{t_0} \frac{t}{A(t)} dt \right)ds.$$
It is unique, up to addition of an affine-linear function of the variables $x$, $xy$.
\end{corollary}

With a similar argument as in the generic case, see Proposition~\ref{coroOrtPolyt=orthoMet}, we derive from the Proposition~\ref{defnCALABIlocalMet} the following proposition.
\begin{proposition}\label{propCalab=trap} Let $(M,\omega, T)$ be a symplectic toric $4$--orbifold with moment map $\mu$. There exists a a compatible Calabi toric metric on $(M,\omega)$ if and only if $\im\mu$ is a trapezoid which is not a parallelogram.
\end{proposition}

\subsection{Extremal Calabi toric metrics}\label{sectCalabiMetric}

Using the Abreu formula~(\ref{ActionAnglemetric}) and the expression of $g=g_{A,B}$ with respect to the moment coordinates $\sigma_1=x$ and $\sigma_2=xy$, one computes that \begin{equation}
  S(G_{g_{A,B}}) = - \sum_{ij}\frac{\del^2}{\del\sigma_i\del\sigma_j} H_{ij} = - \frac{B''(y) + A''(x)}{x}.
\end{equation} In particular, the extremality condition for such a metric may be expressed as conditions on $A$ and $B$. More precisely, since that, on the manifold $\mathring{M}$, the metric $g_{A,B}$ is a smooth Calabi type metric, we can apply~\cite[Proposition 14]{wsd} to obtain:
\begin{proposition}\cite{calabi,wsd}\label{CalabEXTcondWSD}
Let $(M^4,\omega,J,g_{A,B}, T)$ be a Calabi toric orbifold with coordinates $x,y$. The metric $g_{A,B}$ is extremal if and only if $A$ is a polynomial of degree at most $4$, say $A(x)= A_0x^4 + A_1x^3 + A_2 x^2 + A_3 x + A_4$, and $B$ is a polynomial of degree $2$ such that
\begin{equation}\label{calabiExtcond} B''(y)=-2A_2.
\end{equation} Moreover, assuming~(\ref{calabiExtcond}), $g_{A,B}$ has constant curvature if and only if $A_0=0$ and is K\"ahler--Einstein if
\begin{equation}\label{calabiKE} A_1=A_3=0.
\end{equation}\end{proposition}

In particular, if $g_{A,B}$ is extremal, then
\begin{equation}\label{Scalarcalabi} S(G_{g_{A,B}})= -12A_0 x - 6A_1 =-12A_0 \sigma_1 - 6A_1
\end{equation} is equipoised. Thus, we deduce

\begin{corollary}\label{coroCalabiEquipoised} If the Calabi toric metric $g_{A,B}$ is extremal, then $\zeta_{(\Delta,u)}$ is equipoised and the metric $g_{\Sigma} = \frac{dy^2}{B(y)} + B(y)ds^2$ has constant (positive) scalar curvature $2\kappa = - B''(y).$
\end{corollary}
The following is the counterpart of Lemma~\ref{condEXTlemma} in the Calabi toric case.
\begin{lemma}\label{PropCondExtcalabi}
Let $(\Delta,u)$ be a labeled Calabi trapezoid with Calabi parameters $(\alpha_1,\alpha_2,\beta_1,\beta_2,C_{\alpha_1}, C_{\alpha_2}, C_{\beta_1},C_{\beta_2})\in\bR^8$ There exist polynomials $A$, $B$ of respective degree at most $4$ and $2$, satisfying~(\ref{eq:CondCompactCalabi}) and~(\ref{calabiExtcond}) if and only if
\begin{align}\label{ExtremCalabi}
C_{\beta_2} =- C_{\beta_1}.
\end{align} In that case, the polynomials $A$ and $B$ are uniquely determined by the Calabi parameters. $A$ is of degree $3$ if and only if
\begin{align}\label{cscKCalabi}
\frac{1}{(\alpha_2-\alpha_1)^2}\left(\frac{2\alpha_2 + \alpha_1}{C_{\alpha_1}} + \frac{2\alpha_1+ \alpha_2}{C_{\alpha_2}}\right)=-\frac{1}{(\beta_2-\beta_1)C_{\beta_2}},
\end{align} and satisfies, in addition,~(\ref{calabiKE}) if and only if \begin{align}\label{calabiFEparam} \frac{(2\alpha_1 + \alpha_2)\alpha_2}{C_{\alpha_1}} + \frac{(2\alpha_2 + \alpha_1)\alpha_1}{C_{\alpha_2}}=0\end{align}
\end{lemma}
\begin{proof}From the boundary conditions~(\ref{eq:CondCompactCalabi}) and the fact that $B''(y)=-2\kappa$ must be constant (see Proposition~\ref{CalabEXTcondWSD}), we infer that $B(y)= -\kappa(y-\beta_1)(y-\beta_2).$ Hence,
 \begin{equation}\label{CSTscalCurvCalabQutient}B''(y)=-2\kappa =\frac{4}{(\beta_2-\beta_1)C_{\beta_1}}= -\frac{4}{(\beta_2-\beta_1)C_{\beta_2}}\end{equation} is negative. Similarly, the conditions~(\ref{eq:CondCompactCalabi}) and the fact that $A(x)$ must be a polynomial of degree $4$ (see Proposition~\ref{CalabEXTcondWSD}), lead to
$$A(x)= (x-\alpha_1)(x-\alpha_2)Q_A(x)$$
for some degree $2$ polynomial $Q_A(x) = A_0 x^2 + R_1x+R_0$. Thanks to the boundary condition on $A$ we obtain $3$ linear equations for the variables $A_0$, $R_1$, $R_2$.
\begin{equation*}\begin{split}
  A'(\alpha_1)&= (\alpha_1-\alpha_2) Q_A(\alpha_1)= 2 /C_{\alpha_1},\\
  A'(\alpha_2)&= (\alpha_2-\alpha_1) Q_A(\alpha_2)= 2 /C_{\alpha_2},\\
  A_2 &= R_2 -R_1 (\alpha_1+ \alpha_2) + A_0\alpha_1\alpha_2=\kappa.
\end{split}\end{equation*} The solution is
\begin{equation}\label{A0calabi}
\begin{split}
  R_1 =  \frac{1}{(\alpha_2-\alpha_1)^2}&\left(\frac{2}{C_{\alpha_1}} + \frac{2}{C_{\alpha_2}}\right) -A_0 (\alpha_2+\alpha_1)\\
  R_2 =  \frac{-1}{(\alpha_2-\alpha_1)^2}&\left(\frac{2\alpha_2}{C_{\alpha_1}} + \frac{2\alpha_1}{C_{\alpha_2}}\right) + A_0 \alpha_2\alpha_1\\
 A_0 ((\alpha_2+\alpha_1)^2 +2\alpha_1\alpha_2) & = \frac{2}{(\alpha_2-\alpha_1)^2}\left(\frac{2\alpha_2 + \alpha_1}{C_{\alpha_1}} + \frac{2\alpha_1+ \alpha_2}{C_{\alpha_2}}\right)+\kappa.
\end{split}
\end{equation} Notice that $(\alpha_2+\alpha_1)^2 +2\alpha_1\alpha_2>0$ since $\alpha_2>\alpha_1>0$ by assumption.
Thus, we can derive the expression of coefficients $A_0$, $A_1$, $A_2=\kappa$, $A_3$ and $A_4$ in terms of Calabi parameters. The equation~(\ref{cscKCalabi}) is equivalent to $A_0=0$ while~(\ref{calabiKE}) is equivalent to~(\ref{calabiFEparam}).
\end{proof}

\section{Proof of the main results}
\subsection{Equipoised functions and existence of solutions}\label{sectExistSoln=sigma1Parall}

As explained in the introduction, the solution of Problem~\ref{problem1} for parallelograms is the product of solutions on labeled intervals given by formula~(\ref{stadSOLSimplexe}). Notice also that any affine-linear functions are equipoised on a parallelogram. Indeed, this is true on a square and any parallelogram is affinely equivalent to a square.\\

We have established in the previous sections that any quadrilateral which is not a parallelogram is affinely equivalent to either an orthotoric polytope or a Calabi polytope. More precisely, for any convex quadrilateral $\Delta$ which is not a parallelogram there exists an affine invertible map sending $\Delta$ to a quadrilateral in $\bR^2$ which is the image of a rectangle $[\alpha_1,\alpha_2]\times[\beta_1,\beta_2]$ via the map $\sigma(x,y)=(x+y,xy)$ if $\Delta$ is generic and $\sigma(x,y)=(x,xy)$ if $\Delta$ is a trapezoid, see Definitions~\ref{defnOrthoPolyt},~\ref{defnCalabiPolyt}. We fix such a representative of $\Delta$ and a set of $\Vlabel$s $u= (u_1,u_2, u_3, u_4)$. Any other choice of normal inward vectors is obtained from $u$ by
\begin{equation}\label{setofLabel}
  u(r)=\left(\frac{u_1}{r_1},\frac{u_2}{r_2}, \frac{u_3}{r_3}, \frac{u_4}{r_4}\right)
\end{equation}
with $r\in\bR^4_{> 0}$. In the orthotoric (resp. Calabi) case, we choose $u$ such that $u(r)$ give the $\Vlabel$s (\ref{labelOrtho}) (resp. (\ref{labelCalabi})) for
\begin{align}\label{NotationCr}
r= (r_1,r_2,r_3,r_4)=\left(\frac{1}{C_{\alpha_1}},\frac{-1}{C_{\alpha_2}},\frac{-1}{C_{\beta_1}},\frac{1}{C_{\beta_2}}\right).
\end{align}

\subsubsection{The space of normal inward vectors and formal solutions}

In Sections~\ref{sectOrtho} and \ref{sectCalabi}, we introduced two kinds of solution of equation~(\ref{extrem}) in terms of a $S^2\bR^2$--valued functions $\bfH_{A,B}$ depending on two polynomials $A$ and $B$ respectively given by~(\ref{bfHABortho}) if $\Delta$ is generic and~(\ref{bfHABcalabi}) if $\Delta$ is a trapezoid. In both cases, we obtained a criterion on $\Vlabel$s $u(r)$ for $(\Delta,u(r))$ to admit such a solution, expressed as a homogeneous equation with respect to the variables $(r_1,r_2, r_3, r_4)$, see Lemmas~\ref{condEXTlemma} and~\ref{PropCondExtcalabi}. Let $\bfE(\Delta)$ be the cone of inward $\Vlabel$s (labeling $\Delta$) satisfying this criterion. More precisely, using the parametrization above, if $\Delta$ is generic
\begin{equation}
\bfE(\Delta)= \{(r_1,r_2,r_3,r_4) \in \bR_{>0}^4 \,|\; \mbox{(\ref{condEXTpoly5}) holds }\},
\end{equation}
and if $\Delta$ is a trapezoid
\begin{equation}
\bfE(\Delta)= \{(r_1,r_2,r_3,r_4)\in \bR_{>0}^4 \,|\; r_3=r_4\}.
\end{equation}

\begin{lemma}\label{Ext=3dim}
For any convex quadrilateral $\Delta$ which is not a parallelogram, $\bfE(\Delta)$ is a $3$--dimensional cone.
\end{lemma}
\begin{proof} Letting $\bfE'(\Delta)$ be the vector space spanned by $\bfE(\Delta)$, we have $\bfE(\Delta)= \bR_{>0}^4\cap\bfE'(\Delta).$ If $\Delta$ is a trapezoid, the Lemma~\ref{Ext=3dim} obviously follows from the definition and Proposition~\ref{PropCondExtcalabi}. If $\Delta$ is generic, we denote by $(\alpha,\beta)$ its characteristic pair, see Corollary~\ref{coroCharactPair}, and identify $\Delta$ with $\sigma\big([1,\alpha]\times [0,\beta]\big)$. Using the notation~(\ref{NotationCr}), equation~(\ref{condEXTpoly5}) becomes

\begin{equation}\label{condEXTpoly5redD}
D_1\frac{r_1}{(\alpha-1)^2}- D_2\frac{r_2}{(\alpha-1)^2} -D_3\frac{r_3}{\beta^2} + D_4\frac{r_4}{\beta^2} = 0,
\end{equation}where
\begin{equation*}
\label{condEXTpoly5red}
\begin{split}
D_1 &= (1 +\alpha)^2+2\alpha^2 +\beta^2 -2(2\alpha + 1)(\beta), \\
D_2 &= (1 +\alpha)^2+2 +\beta^2 -2(2 + \alpha)\beta,\\
D_3 &= 3\beta^2 +(1 +\alpha)^2 +2\alpha -4\beta(1 +\alpha),\\
D_4 &= \beta^2 + (1 +\alpha)^2 +2\alpha -2\beta(1+\alpha).
\end{split}
\end{equation*}The lemma follows from the fact that $D_1>D_4>D_3>D_2>0$ which in turn follows that $0<\beta<1<\alpha$ as stated in Corollary~\ref{coroCharactPair}.
\end{proof}

\subsubsection{The space of normal inward vectors with equipoised extremal affine function}

Recall that an affine function $f$ is \defin{equipoised} on $\Delta$ if $\sum_{i=1}^4 (-1)^i f(s_i) =0$, where $s_1,\dots,s_4$ are the vertices of $\Delta$ and $s_1$ is not adjacent to $s_3$. Moreover, we can assume (without loss of generality, see Lemmas~\ref{theoNoParall=ortho}~\ref{lemTrap=Calab}) that any quadrilateral, unless parallelogram, is either an orthotoric polytope or a Calabi polytope. Note that orthotoric and Calabi polytopes come with a fixed affine embedding $\Delta\subset\bR^2=\{(\sigma_1,\sigma_2)\,|\;\sigma_i \in \bR\}$ having as common feature that the function $(\sigma_1,\sigma_2)\mapsto\sigma_1$ is equipoised on $\Delta$ while $(\sigma_1,\sigma_2)\mapsto\sigma_2$ is not. Hence, an affine-linear function $f\co\Delta \ra \bR$ is equipoised on $\Delta$ if and only if $f$ is constant with respect to $\sigma_2$.

Hence, from the linear system~(\ref{systCHPext}), for a given labeling $u= (u_1,u_2, u_3, u_4)$, the extremal affine function $\zeta_{(\Delta,u)}$ is equipoised if and only if the linear system
\begin{equation}\label{systCHPextParall}
\begin{split}
W_{01} \zeta_1  + W_{00} \zeta_0 &= Z_0\\
W_{11} \zeta_1  + W_{01} \zeta_0 &= Z_1\\
W_{21} \zeta_1  + W_{02} \zeta_0 &= Z_2
\end{split}
\end{equation}
admits a solution, that is if and only if
\begin{equation}\label{condChpExtParall}
( W_{00}W_{11}-W_{01}^2)Z_2 = (W_{02}W_{11}-W_{12}W_{01})Z_0  +(W_{12}W_{00}-W_{02}W_{01})Z_1.
\end{equation}

Notice that in the system~(\ref{systCHPext}) (and hence~(\ref{systCHPextParall})), only the right-hand side depends on the $\Vlabel$s: For any other inward vectors $u(r)$, see~(\ref{setofLabel}), with $r\in\bR^4_{> 0}$, the matrix $W$ remains unchanged while $Z(r)$ depends linearly on $r=(r_1,r_2,r_3,r_4)$. Hence, the equation~(\ref{condChpExtParall}) is homogeneous and linear with respect to $r$. Denote $a_{ij} = \int_{F_j} \mu_i d\nu_j$ where $u_j \wedge d\nu_j= d\mu_1\wedge d\mu_2$ on $F_j$, so that $Z_i(r) = \sum_j a_{ij}r_j$. The equation~(\ref{condChpExtParall}) can be rewritten as $\sum_{j=1}^4 E_j r_j=0$, where
$$E_j=a_{2j} ( W_{00}W_{11}-W_{01}^2) -  a_{0j}(W_{02}W_{11}-W_{12}W_{01}) -a_{1j}(W_{12}W_{00}-W_{02}W_{01}).$$
We define the vector space $\bfB'(\Delta)= \left\{ r\in\bR^4 \left|\; \sum_{j=1}^4 E_j r_j=0\right.\right\},$
so that $$\bfB(\Delta)=\bfB'(\Delta)\cap \bR_{>0}^4= \{ r\in\bR_{>0}^4\,|\; \zeta_{(\Delta,u(r))} \mbox{ is equipoised on }\Delta\}.$$ The next Lemma will be the key for the main result of this paper. It follows from the fact that we linearized the PDE of Problem~\ref{problem1} in Lemmas~\ref{condEXTlemma} and~\ref{PropCondExtcalabi}.

\begin{lemma}\label{E=B}
For any convex quadrilateral $\Delta$ which is not a parallelogram $$\bfB(\Delta) =\bfE(\Delta).$$
\end{lemma}

\begin{proof} Since $\bfB(\Delta)=\bfB'(\Delta)\cap \bR_{>0}^4$, $\bfE(\Delta)=\bfE'(\Delta)\cap \bR_{>0}^4$ and since we proved that $\bfE'(\Delta)$ is a $3$--dimensional vector space, we only have to prove that $\bfE(\Delta)\subset \bfB(\Delta)$ and that $\bfB'(\Delta)$ is a $3$--dimensional vector space.

Firstly, $\bfE(\Delta)\subset \bfB(\Delta)$. Indeed, for a point $r\in\bfE(\Delta)$ one can construct a matrix-valued function $\bfH_{A,B}=(H_{ij})$ satisfying the boundary conditions~(\ref{condCOMPACTIFonH}) and~(\ref{extrem}) by Lemmas~\ref{condEXTlemma} and~\ref{PropCondExtcalabi}. The latter condition means that the function
\begin{equation}\label{curvH}
  S(\bfH_{A,B})=-\frac{\del^2H_{ij}}{\del\mu_i\del\mu_j}
\end{equation} is affine-linear and thus coincides with $\zeta_{(\Delta,u(r))}$, see Remark~\ref{dependenceExtchp2}. In particular, due to the construction $\zeta_{(\Delta,u(r))}$ is equipoised and then $r\in \bfB(\Delta)$, see Remark~\ref{sigma1parall} and Corollary~\ref{coroCalabiEquipoised}.

Secondly, $\bfB'(\Delta)$ has dimension $3$ since the coefficients $E_i$ of the defining equation $\sum_{j=1}^4 E_j r_j=0$ are not all zero. We prove this fact separately for both cases.

 $(i)$ Suppose that $\Delta$ is generic, we identify $\Delta$ with an orthotoric polytope with parameters $\beta_1<\beta_2<\alpha_1<\alpha_2$. Recall that the (fixed) inward normal vectors $u$ are given via~(\ref{setofLabel}) and~(\ref{NotationCr}), that is,
$$u_1=\begin{pmatrix}
        \alpha_1\\
        -1
      \end{pmatrix},\;\; u_2=\begin{pmatrix}
        \alpha_2\\
        -1
      \end{pmatrix},\;\; u_3=\begin{pmatrix}
        \beta_2\\
        -1
      \end{pmatrix},\;\; u_4=\begin{pmatrix}
        \beta_2\\
        -1
      \end{pmatrix}.$$
We compute
\begin{equation} \label{expressParamZ}
\begin{split}
 Z_0(r)&= 2(\beta_2-\beta_1) \left( r_2 + r_1 \right) + 2(\alpha_2 - \alpha_1) \left( r_4+r_3 \right) \\
Z_1(r) &= (\beta_2 -\beta_1)\left((2\alpha_2 +\beta_2 +\beta_1)r_2 + (2\alpha_1 +\beta_2 +\beta_1)r_1 \right) \\
&\quad  +(\alpha_2 -\alpha_1)\left((2\beta_2 +\alpha_2 +\alpha_1)r_4 + (2\beta_1 +\alpha_2 +\alpha_1)r_3 \right)\\
Z_2(r) &= (\beta^2_2 -\beta^2_1)\left(\alpha_2r_2 + \alpha_1r_1 \right) + (\alpha^2_2 -\alpha^2_1)\left(\beta_2r_4 + \beta_1r_3 \right).
\end{split}
\end{equation}
In particular, at least one coefficient is not zero as the alternating sum $-E_1 +E_2 -E_3 +E_4$ is
\begin{equation}\label{coeffNonTousNuls}
(\alpha_2-\alpha_1)(\beta_2-\beta_1)(W_{00}W_{11}-W_{01}^2)\left(\alpha_2 +\alpha_1 -\beta_2 -\beta_1\right) \neq 0.
\end{equation}
Indeed, recall that $W_{ij} =\int_{\Delta}\mu_i\mu_j dv$ is an inner product on $L^2(\Delta)$ and that $\mu_0$ and $\mu_1$ are everywhere independent as functions on $\bR^2$. Hence, the Cauchy--Schwarz inequality is strict for them, showing that $W_{01}^2>W_{00}W_{11}$. Thus, (\ref{coeffNonTousNuls}) holds since $\beta_1<\beta_2<\alpha_1<\alpha_2$.

 $(ii)$ Suppose that $\Delta$ is a trapezoid, it corresponds to a Calabi polytope with parameters $\beta_2>\beta_1\geq 0$ and $\alpha_2>\alpha_1>0$. The fixed inward normals are
$$u_1=\begin{pmatrix}
        \alpha_1\\
        0
      \end{pmatrix},\;\; u_2=\begin{pmatrix}
        \alpha_2\\
        0
      \end{pmatrix},\;\; u_3=\begin{pmatrix}
        \beta_2\\
        -1
      \end{pmatrix},\;\; u_4=\begin{pmatrix}
        \beta_2\\
        -1
      \end{pmatrix}$$
and we compute
\begin{align*}  Z_0(r)=& (\beta_2-\beta_1) \left(2r_2 +2r_1 \right) + (\alpha_2^2 - \alpha_1^2) \left(r_4 + r_3 \right) \\
Z_1(r) =& (\beta_2 -\beta_1)\left(2\alpha_2r_2 + 2\alpha_1r_1 \right) + \frac{2}{3}(\alpha_2^3 - \alpha_1^3) \left(r_4+ r_3 \right)\\
Z_2(r) =& (\beta^2_2 -\beta^2_1)\left(\alpha_2r_2 +\alpha_1r_1 \right) + \frac{2}{3}(\alpha^3_2 -\alpha^3_1)\left(\beta_2r_4 + \beta_1r_3 \right).
\end{align*}
Notice that
$$E_4-E_3= -\frac{2}{3}(\beta_2-\beta_1)(W_{00}W_{11}-W_{01}^2)(\alpha_2^3 - \alpha_1^3)  \neq 0$$
as above. It follows that at least one coefficient is non-zero.
\end{proof}

\subsection{Proof of Theorem~\ref{mainTHEOparall}}\label{sectStabilOrtho}

The proof given below is similar to the one of~\cite[Theorem 3]{ambi}.

Let $(\Delta, u)$ be a labeled convex quadrilateral with a equipoised $\zeta_{(\Delta,u)}$. Suppose that $\Delta$ is not a parallelogram. From Lemma~\ref{E=B} we know that there exist polynomials $A$ and $B$ for which the matrix $\bfH_{A,B}$, given by~(\ref{bfHABcalabi}) for trapezoids and~(\ref{bfHABortho}) for generic quadrilaterals, satisfies the boundary conditions~(\ref{condCOMPACTIFonH}) and $S(\bfH_{A,B})=\zeta_{(\Delta,u)}$. Thanks to the uniqueness of the solution of Problem~\ref{problem1}, see~\S\ref{uniqueness}, to prove Theorem~\ref{mainTHEOparall} we only have to show that {\it $\bfH_{A,B}$ is positive definite if and only if $(\Delta,u)$ is relatively analytically $K$--stable with respect to toric degenerations.}\\

Denote by $\mP(\Delta)$ the set of strictly convex piecewise affine-linear functions on $\Delta$. Suppose that $f\in\mP(\Delta)$ has only one crease, meaning that $\Delta$ is cut into two pieces $\Delta=\Delta_1\cup\Delta_2$ on which $f$ is affine-linear (thus $f_i=f_{|\Delta_i}$ is smooth). Denote by $S_f=\del\Delta_1\cap \del\Delta_2$, the segment passing through the interior of $\Delta$, determined by the crease of $f$. A normal vector to $S_f$ is given by
$$u_f =
\begin{pmatrix}
\frac{\del}{\del \mu_1}f_1 -\frac{\del}{\del \mu_1}f_2\\
 \frac{\del}{\del \mu_2}f_1-\frac{\del}{\del \mu_2}f_2
\end{pmatrix} =\begin{pmatrix} a_{11} -a_{21} \\
a_{12} -a_{22}\end{pmatrix}$$
where $a_{sj} = \frac{\del}{\del\mu_j} f_s$ are constants. This vector is inward to $\Delta_1$. Denote by $d\nu_f$ the (positive) volume form on the oriented segment $S_f$ for which $u_f\wedge d\nu_f = dv=d\mu_1\wedge d\mu_2$.
Substituting $\zeta_{(\Delta,u)}$ by $S(\bfH_{A,B})$ (even though $\bfH_{A,B}$ is not necessarily positive definite), we get
\begin{equation}
\begin{split}\label{creaseFutaki}
\mL_{(\Delta,u)}(f) = \int_{S_f} \bfH_{A,B}(u_f,u_f)\,d\nu_f.
\end{split}
\end{equation}

 For functions with more than one crease, the relative Futaki functional decomposes into a sum over the creases of expressions of the type~(\ref{creaseFutaki}). Indeed, the integration by part leading to~(\ref{creaseFutaki}) may be used successively, as in~\cite{zz}, showing that if $\bfH_{A,B}$ is positive definite then $(\Delta,u)$ is relatively analytically $K$--stable with respect to toric degenerations. It remains to prove the converse.\\

From Lemmas~\ref{theoNoParall=ortho} and~\ref{lemTrap=Calab}, we can assume that $\Delta$ is the image of a rectangle $[\alpha_1,\alpha_2]\times[\beta_1,\beta_2]$ via the affine map $\sigma$, see Definitions~\ref{defnOrthoPolyt},~\ref{defnCalabiPolyt}. Let $C_{\alpha_1}$, $C_{\alpha_2}$, $C_{\beta_1}$ and $C_{\beta_2}$, be the constant determining the $\Vlabel$s respectively as in~(\ref{labelOrtho}) (resp. (\ref{labelCalabi})). Then, to every $x \in [\alpha_1,\alpha_2]$ corresponds a segment $S_x\subset\Delta$, given by the image of $[\beta_1,\beta_2]$ by $\sigma(x,\cdot)$. We define similarly segments $S_y$ for every $y\in [\beta_1,\beta_2]$. For any function $f$ whose only crease is $S_x$, we compute
\begin{equation}\label{distcreaseA}
 \mL_{(\Delta,u)}(f) =\int_{S_{x}}\bfH_{A,B}(u_f,u_f)\,d\nu_{f} = A(x)(\beta_2-\beta_1),
\end{equation}
up to a positive multiplicative constant. Similarly, if $S_y$ is the only crease of $f$, we get that, up to a positive multiplicative constant,
\begin{equation}\label{distcreaseAcalab}
  \mL_{(\Delta,u)}(f) =\int_{S_{y}}\bfH_{A,B}(u_f,u_f)\,d\nu_{f} = B(y) (\alpha_2-\alpha_1).
\end{equation}

Then, from~(\ref{distcreaseA}) and~(\ref{distcreaseAcalab}), $\mL_{(\Delta,u)}(f) >0$ for any non-affine-linear function $f\in\mP(\Delta)$ implies that $A$ and $B$ are positive on the respective open intervals $(\alpha_1,\alpha_2)$ and $(\beta_1,\beta_2)$, and thus $\bfH_{A,B}$ is positive definite.

\begin{rem} Of special interest is the case when the parameters $\alpha_1$, $\alpha_2$, $\beta_1$, $\beta_2$, $C_{\alpha_1}$, $C_{\alpha_2}$, $C_{\beta_1}$, $C_{\beta_2}$ are rational. One then gets a rational labeled polytope $(\Delta,\Lambda, u)$ with vertices lying in $\Lambda^*$, see Lemma~\ref{lemstronglyRatIFFrat}. Presumably, this would be the case where the algebro-geometric setting of the problem makes sense, see~\cite{don:scalar} for the case of smooth varieties and the work of~\cite{RTorbi} for the case of orbifolds with cyclic orbifold structure groups. In this setting, \defin{rational} convex piecewise functions would arise from toric degenerations of the toric orbifold.

It is worth noticing, as in~\cite[Theorem 3]{ambi}, that under this assumption the polynomials $A$ and $B$ have rational coefficients and two simple rational roots. Therefore, they cannot admit double irrational roots, showing that $H_{A,B}$ is definite positive as soon as $\mL_{(\Delta,u)}(f)>0$ for any rational function in $\mP(\Delta)$. \end{rem}

\subsection{Proof of Theorem~\ref{propCONEsol}}

We prove Theorem~\ref{propCONEsol} in three steps. \newline
First, we show the last part of the statement concerning $\bfC(\Delta)$ and $\bfK(\Delta)$. Then we prove that $\bfB^+(\Delta)$ and $\bfB(\Delta)\backslash \bfB^+(\Delta)$ are both non-empty for generic quadrilaterals and, finally, $\bfB^+(\Delta)$ is non-empty for trapezoids. We begin by recalling and giving alternative definitions for $\bfB^+(\Delta)$, $\bfC(\Delta)$ and $\bfK(\Delta)$.

Let $\Delta$ be a convex quadrilateral which is not a parallelogram. We suppose, without loss of generality, that $\Delta$ is embedded in $\bR^2$ as an orthotoric quadrilateral if $\Delta$ is generic and as a Calabi trapezoid otherwise. In both cases the associated parameters are denoted $\alpha_1$, $\alpha_2$, $\beta_1$, $\beta_2$.

Let $\bfN(\Delta)$ be the $4$--dimensional cone of inward normals associated to the facets of $\Delta$. We fix inward normal vectors $u= (u_1,u_2, u_3, u_4)$ so that any other normal inward vectors can be expressed as $u(r)$ for some $r\in\bR^4_{> 0}$, using conventions (\ref{setofLabel}) and~(\ref{NotationCr}). Via this parametrization, $\bfN(\Delta)= \bR_{>0}^4$.\\

We know from Lemmas~\ref{Ext=3dim} and~\ref{E=B} that the condition that $\zeta_{(\Delta,u)}$ is equipoised defines a codimension one sub-cone $\bfB(\Delta)\subset \bfN(\Delta)$ which can be equivalently defined as the cone $\bfE(\Delta)$ of $\Vlabel$s for which there exist polynomials of degree $4$, $A$ and $B$, such that $\bfH_{A,B}$ (given by~(\ref{bfHABortho}) if $\Delta$ is generic and~(\ref{bfHABcalabi}) otherwise) is a solution of~(\ref{extrem}) satisfying the compactification condition~(\ref{condCOMPACTIFonH}). Recall that $A$ and $B$ are uniquely defined by these conditions. Letting $A(x) =A_0x^4 +A_1x^3+ A_2x^2 + A_3x +A_4$ and going back to the proof of Lemma~\ref{E=B}, we get  \begin{equation*}
  \bfB(\Delta) = \left\{r\in \bfN(\Delta) \left| \zeta_{(\Delta,u(r))}=-12A_0 \sigma_1 -6 A_1\right.\right\}.
\end{equation*} In particular, noticing that a constant function is equipoised on any quadrilateral, the set of $\Vlabel$s $\bfC(\Delta)$ for which $\zeta_{(\Delta,u)}$ is constant is a subset of $\bfB(\Delta)$, and is equivalently defined as \begin{equation*} \bfC(\Delta) = \left\{r\in \bfB(\Delta)\, \left|\; A_0=0 \right.\right\}.\end{equation*} The subset $\bfB^+(\Delta)$ of normals for which $(\Delta,u)$ is relatively analytically $K$--stable corresponds, via Theorem~\ref{mainTHEOparall}, to the subset of $\bfB(\Delta)$ for which $A$ and $B$ are positive respectively on $(\alpha_1,\alpha_2)$, $(\beta_1,\beta_2)$.\\

\noindent {\bf Step 1.} For any $\Vlabel$s in $\bfC(\Delta)\subset \bfB(\Delta)$, $A$ is of degree $3$ and has $\alpha_1$, $\alpha_2$ as roots. The conditions $A'(\alpha_1)>0$ and $A'(\alpha_2)<0$ ensure that $A$ is positive on $(\alpha_1,\alpha_2)$. For similar reasons $B$ is positive on $(\beta_1,\beta_2)$. Hence, $\bfC(\Delta)\subset \bfB^+(\Delta)$. \\

By using Lemmas~\ref{condEXTlemma} and~\ref{PropCondExtcalabi} the sets $\bfC(\Delta)$, $\bfB(\Delta)$ are defined by linear equations with respect to $r_1,r_2,r_3,r_4$. We already know from Lemma~\ref{Ext=3dim} that $\bfB(\Delta)$ is $3$--dimensional. Similarly, the set of $\Vlabel$s $\bfC(\Delta)$ for which $\zeta_{(\Delta,u)}$ is constant is defined by\begin{equation*}
  \bfC(\Delta) = \left\{r\in \bfN(\Delta) \left| \begin{array}{l}
    \mbox{(\ref{condEXTpoly5}), (\ref{anciennement-ii}) hold, if $\Delta$ is generic}\\
    \mbox{(\ref{ExtremCalabi}), (\ref{cscKCalabi}) hold, if $\Delta$ is a trapezoid}
  \end{array} \right.\right\},
\end{equation*} while the set of K\"ahler--Einstein metrics is
 \begin{equation*}
  \bfK(\Delta) = \left\{r\in \bfN(\Delta) \left| \begin{array}{l}
    \mbox{(\ref{condEXTpoly5}), (\ref{anciennement-ii}), (\ref{anciennement-iv}) hold, if $\Delta$ is generic}\\
    \mbox{(\ref{ExtremCalabi}), (\ref{cscKCalabi}), (\ref{calabiFEparam}) hold, if $\Delta$ is a trapezoid}
  \end{array} \right.\right\}.
\end{equation*}
We infer from Lemma~\ref{E=B} that $\bfB(\Delta)$, $\bfC(\Delta)$ and $\bfK(\Delta)$ are sub-cones of $\bfN(\Delta)$ of respective codimension one, two and three.\\

It remains to prove that $\bfB^+(\Delta)$ is proper and is a non-empty open subset of $\bfB(\Delta)$. It is clearly open by definition. Recall from the proofs of Lemmas~\ref{condEXTlemma},\ref{PropCondExtcalabi} that $A(x) = (x-\alpha_1)(x-\alpha_2)Q_A(x)$ and $B(y)= (y-\beta_1)(y-\beta_2)Q_B(y) $ where $Q_A(x) = A_0x^2+R_1x +R_2$ and the degree of $Q_B$ depends whether or not $\Delta$ is generic. The polynomial $A$ (resp. $B$) is positive on $(\alpha_1, \alpha_2)$ (resp. $(\beta_1, \beta_2)$) if and only if $Q_A$ (resp. $Q_B$) is negative on these intervals. The compactification conditions imply that $Q_A$ (resp. $Q_B$) is negative at the ends of the interval $(\alpha_1, \alpha_2)$ (resp. $(\beta_1, \beta_2)$). In particular, if $A_0>0$ then $A$ is positive on $(\alpha_1, \alpha_2)$.\\

\noindent {\bf Step 2.} If $\Delta$ is generic then $Q_B(y)= -A_0y^2 +S_1y+S_2$. The fact that $\bfB^+(\Delta)$ is a non-empty open subset of $\bfB(\Delta)$ will follow if we can find $r\in\bfB(\Delta)$ for which $A_0>0$ and $Q_B$ has imaginary roots. Indeed, $A_0>0$ implies $A>0$ on $(\alpha_1,\alpha_2)$ as above, and since $Q_B$ has no real root and is negative at $\beta_1$, $Q_B$ is always negative. Thus, $r$ is in the open subset of $\bfB(\Delta)$ (included in $\bfB^+(\Delta)$) defined by $A_0>0$ and $S_1^2 +4S_2A_0<0$. On the other hand, the fact that $\bfB^+(\Delta)$ is a proper subset of $\bfB(\Delta)$ would follow from the existence of $r \in\bfB(\Delta)$ for which $A_0>0$ and $Q_B$ has a double root in $(\beta_1, \beta_2)$. We now show the existence of such $r$.\\

Assume, (without loss of generality, see Corollary~\ref{coroCharactPair}), that $\Delta$ is the orthotoric quadrilateral with characteristic pair $(\alpha,\beta)$, where $0<\beta<1<\alpha$ and $\alpha-\beta \geq 1$. We use the notation in the proof of Lemma~\ref{Ext=3dim}. Let $r=(r_1,r_2,r_3,r_4)\in \bR_{>0}^4$. Letting $a=\frac{r_2-r_1}{(\alpha-1)^2}$ and $b=\frac{r_4-r_3}{\beta^2}$, $r\in \bfB(\Delta)$ if and only if \begin{equation}
\label{v1}\frac{r_1}{(\alpha-1)^2}(D_1-D_2) = -\frac{r_3}{\beta^2}(D_4-D_3) +aD_2 -bD_4.
\end{equation} Recall from Lemma~\ref{Ext=3dim} that $D_1>D_4>D_3>D_2>0.$\\

Now $r\in \bfB^+(\Delta)$ and $A_0>0$ if and only if, assuming~(\ref{v1}), the following conditions hold:
 \begin{align}
   & a>-\frac{r_1}{(\alpha-1)^2}, \; \frac{r_3}{\beta^2}>\mbox{max}\{0,-b\} \;\mbox{ and }\; \frac{r_3}{\beta^2}(D_4-D_3) < aD_2 -bD_4, \label{anc1}\\
   & a>b. \label{anc2}
 \end{align} Moreover, $Q_B$ has conjugate imaginary roots if and only if $S_1^2 + 4A_0S_2 \leq0$, that is, if and only if
 \begin{equation} \label{anc3} \frac{r_3}{\beta^2} \geq \frac{b^2(\alpha+1-\beta)}{2(a-b)\beta} +\frac{(a-b)\beta}{8(\alpha+1-\beta)} -b,\end{equation} with equality if and only if $Q_B$ has the double (real) root $\lambda=\frac{S_1}{2A_0}$.

\begin{claim}\label{claim} For any $b<0$, there exists $a_0>0$ such that for all $a>a_0$, we have \begin{equation} \label{IneGv3}\frac{1}{D_4-D_3}\left(aD_2 -bD_4\right)>\frac{b^2(\alpha+1-\beta)}{2(a-b)\beta} +\frac{(a-b)\beta}{8(\alpha+1-\beta)} -b.\end{equation}
 \end{claim}

  Step 2 will be complete as soon as we prove the Claim~\ref{claim}. Indeed, taking $b$ and $a>a_0$ satisfying this Claim we have that $a>0>b$ implies (\ref{anc2}) and the inequality~(\ref{IneGv3}) allows us to pick $r_3$ such that \begin{align}&\frac{r_3}{\beta^2} <\frac{1}{D_4-D_3}\left(aD_2 -bD_4\right), \label{IneG+}\\
 &\frac{r_3}{\beta^2}\geq\frac{b^2(\alpha+1-\beta)}{2(a-b)\beta} +\frac{(a-b)\beta}{8(\alpha+1-\beta)} -b.\label{IneG+2}\end{align} Thus, (\ref{IneG+2}) implies that $\frac{r_3}{\beta^2}>-b$ while (\ref{IneG+}) together with equation~(\ref{v1}) imply that $r_1>0$, so the condition (\ref{anc1}) holds.
 Picking $r_3$ such that the inequality (\ref{IneG+2}) is strict implies that the strict inequality of (\ref{anc3}) holds (so that $Q_B$ has no real roots) while picking $r_3$ such that the equality of (\ref{IneG+2}) holds implies that $Q_B$ has the double root $\lambda$. We have $\lambda \in (0,\beta)$ if and only if $S_1 \in (0,2A_0\beta)$ which means (by virtue of~(\ref{definETA})) that $-A_0\beta < 2b <A_0 \beta$. Thus, $Q_B$ has a double root in $(0,\beta)$ if and only if \begin{equation}\label{Inegdoubleroot}
    -\frac{\beta(a-b)}{2(\alpha+1-\beta)} < 2b < \frac{\beta(a-b)}{2(\alpha+1-\beta)}
\end{equation} which, in turn, is verified as soon as $a$ is big enough. Hence, it remains to prove the Claim~\ref{claim}.
\begin{rem} It is easy to see that once given $r\in\bfB(\Delta)\backslash \bfB^+(\Delta)$ such that $A_0>0$ and $Q_B$ has a double root, there is a $3$--parameters family of such solutions.
\end{rem}
\begin{proof}[Proof of the Claim~\ref{claim}] We compute that $$\frac{1}{D_4-D_3}\left(aD_2 -bD_4\right) = \frac{(a-b)[(\alpha-\beta)^2 +3 +2\alpha-4\beta] +2b(1-\beta-\alpha)}{2\beta(\alpha+1-\beta)}.$$
One can assume that $a>b$ so the inequality~(\ref{IneGv3}) holds if and only if \begin{equation}\begin{split}&\left[\frac{4(\alpha-\beta)^2 +8\alpha+12-16\beta -\beta^2 }{8\beta(\alpha+1-\beta)} \right](a-b)^2\\
 &+ \left[\frac{2b(1-\beta-\alpha)}{2\beta(\alpha+1-\beta)} +b\right](a-b) -\frac{b^2(\alpha+1-\beta)}{2\beta}>0.\end{split}\end{equation} Fixing $b<0$ the left hand side is a polynomial, say $P(a-b)$, of degree two with respect to $a-b$ for which the main coefficient $$\frac{4(\alpha-\beta)^2 +8\alpha+12-16\beta -\beta^2 }{8\beta(\alpha+1-\beta)}$$ is positive. Thus $P$ is a convex function. Hence, there exists $a>0$ big enough to ensure $P$ to be positive at $a-b$.
\end{proof}

\noindent {\bf Step 3.} If $\Delta$ is a trapezoid, then $Q_B(y)= -\kappa= -2r_3$ and $\bfB(\Delta)= \{r\in\bR^4_{>0}\,|\; r_3=r_4\}$. Using the formulae of Lemma~\ref{PropCondExtcalabi}, we express $A_0$ in terms of the variables $r=(r_1,r_2,r_3,r_4)$ and the Calabi parameters $\alpha_1$, $\alpha_2$, $\beta_1$, $\beta_2$. In particular, for any number $r_3=r_4>0$, there exist $r_1,r_2>0$ such that $(r_1,r_2,r_3,r_4)\in \bfB(\Delta)$ and $A_0 >0$. Thus, for such $r$, $A$ is positive on $(\alpha_1,\alpha_2)$. We then infer that $\bfB^+(\Delta)$ is not empty.

\begin{rem} \label{remWBF}The classification presented in this paper provides automatically a classification of toric weakly Bochner-flat metrics (i.e with co-closed Bochner tensor). Indeed, weakly Bochner-flat metrics are extremal and an alternative definition is that $(g, J,\omega)$ is weakly Bochner-flat if the form $$\tilde{\rho}_g = \rho_g - \frac{Scal_g}{2m(m+1)}\omega$$ where $\rho_g$ is the Ricci form, is a Hamiltonian $2$--form, see~\cite{H2FI}. In particular, if $(g, J,\omega)$ is a toric weakly Bochner-flat metric then it is a toric K\"ahler-Einstein metric or admits a non trivial Hamiltonian $2$--form, $\tilde{\rho}_g$.

In the setting of toric geometry, this latter case implies that the moment polytope $\Delta$ is either a triangle or a quadrilateral see~\S\ref{sectGEN=ortho} and~\ref{sectCalabiPolyt}. If $\Delta$ is a triangle, using the uniqueness of extremal metrics the metric $g$ should be a Bochner-flat metric on a weighted projective space as classified in~\cite{bryant}, with symplectic potential given by~(\ref{stadSOLSimplexe}). If $\Delta$ is a quadrilateral, depending of the number of its parallel edges, the metric $g$ is either a product of metrics, a Calabi-type metric or an orthotoric metric. Moreover, if $\Delta$ is a quadrilateral which is not a parallelogram, using again the local characterization of metrics admitting Hamiltonian $2$--form of~\cite{H2FI} the condition of being weakly Bochner-flat metric is a linear condition on the coefficients of the polynomials $A$ and $B$. More precisely, the normals $u$ leading to a (formal) weakly Bochner-flat metric $\bfH_{A,B}$ form a sub-cone in $\bfB(\Delta)$, defined via the linear equation $A_3=-B_3$ if $\Delta$ is generic and $A_3=0$ if $\Delta$ is a trapezoid (without assuming $A_0=0$).
\end{rem}

\section{Geometric applications}\label{sectApplication}

  A labeled polytope $(\Delta, u)$ is associated to a symplectic toric orbifold via the Delzant--Lerman--Tolman correspondence if and only if $(\Delta, u)$ is a rational labeled polytope with respect to a lattice $\Lambda$. The first part of this section gives an intrinsic criterion for testing rationality of polygons. In this paper, {\it polygon} refers to $2$--dimensional polytopes. In particular, they are compact and convex.

\subsection{The rational type condition}\label{classifQuad}\label{sectRationality} \label{sectCrossRatio} \label{sectRationalortho}

Recall that $(\Delta,u)$ is {\it rational with respect to a lattice} $\Lambda$ if $u_i \in \Lambda$ and $\Delta$ is of {\it rational type} if there exists a lattice $\Lambda$ and a set of $\Vlabel$s $u$ such that $(\Delta,u)$ is rational with respect to a lattice $\Lambda$, see Definition~\ref{defnRatpolyt}. Let $\Delta$ be a polytope with $d$ facets in a $2$--dimensional affine space, $(\mA,V)$. There is a canonical way to associate $d$ (not necessarily distinct) points of $\bP(V^*)$: To each facet, we associate its normal line.

Recall that the \defin{cross-ratio} is defined on ordered sets of four distinct points of the real projective line, $P_i = [x_i:y_i] \in \brp^1$, $1\leq i\leq 4$, with $P_i\neq P_j$ if $i\neq j$, by the formula
$$\cratio(P_1,P_2;P_3,P_4) = \frac{(x_1y_3-y_1x_3)(x_2y_4-y_2x_4)}{(x_1y_4-y_1x_4)(x_2y_3-y_2x_3)}.$$
This definition does not depend on the chosen representatives and is invariant under projective transform. The cross-ratio may alternatively be defined for an ordered set of four, non-zero, distinct vectors.
\begin{rem}\label{orderCrossRatio}
For any permutation $\gamma \in S_4$ and four distinct points $P_1$, $P_2$, $P_3$, $P_4\in\bP(V)$, the number $\cratio_{\gamma}=\cratio(P_{\gamma(1)},P_{\gamma(2)};P_{\gamma(3)},P_{\gamma(4)})$ lies in the set
$$\{\cratio, \cratio^{-1}, 1-\cratio, (1-\cratio)^{-1},\frac{\cratio}{\cratio-1},\frac{\cratio-1}{\cratio}\}$$
where $\cratio =\cratio_{\mathrm{id}}=\cratio(P_1,P_2;P_3,P_4)$. Hence, $\cratio_{\gamma}$ is rational if and only if $\cratio$ is.
\end{rem}

\begin{rem} \label{3transitive}
Since $PGL(2,\bR)$ acts simply $3$--transitively on $\brp^1$, for any ordered distinct three points $P_1$, $P_2$, $P_3 \in \brp^1$, there exists a unique $A\in PGL(2,\bR)$ such that $AP_1=[1:0]$, $AP_2=[0:1]$, $AP_3 =[1:1]$. Then, for any $P_4\in\brp^1$
$$\cratio(P_1,P_2;P_3,P_4)=\cratio([1:0],[0:1];[1:1],AP_4)= \mbox{slope} (AP_4)$$
\end{rem}
\begin{proposition}\label{crossratioTHEO} Let $\Delta$ be a polytope with $d$ edges in a $2$--dimensional affine space. $\Delta$ is of rational type if and only if $\Delta$ has either
\begin{enumerate}
\item at most $3$ distinct normal lines,
\item $4$ distinct normal lines with rational cross-ratio,
\item at least $4$ distinct normal lines and the cross-ratio of any four of them is rational.
\end{enumerate}
\end{proposition}

\begin{proof} Let $\Delta$ be a convex polygon in $\bR^2$. For at most $3$ points of $\brp^1$, there obviously exists a lattice intersecting non-trivially each of them, see Remark~\ref{3transitive}. Thus we suppose that $\Delta$ has at least $4$ normal lines $\delta_1$, ... $\delta_k\in\brp^1$.

Suppose there exists a lattice $\Lambda$ intersecting non-trivially $\delta_1$, ... $\delta_k$. Then, there exists $A\in GL(2,\bR)$ such that $A(\Lambda)=\bZ^2$. Since $(A\delta_i \cap \bZ^2)\neq \{0\}$ for all $i$, we can choose a non-zero integral point in each real line $\delta_i$ to compute the cross-ratio: For any four distinct indices $i$, $j$, $k$, $l$, we get
$$\cratio(\delta_i,\delta_j;\delta_k,\delta_l)=\cratio(A\delta_i,A\delta_j;A\delta_k,A\delta_l)\in\bQ.$$

Conversely, fix three of the normal lines, say $\delta_1$, $\delta_2$ and $\delta_3$. By Remark~\ref{3transitive}, there exists a unique $[A]\in PGL(2,\bR)$ such that $[A]\delta_1=[0:1]$, $[A]\delta_2=[1:0]$ and $[A]\delta_3=[1:1]$. Thus, by assumption and Remark~\ref{3transitive},
$$\cratio(\delta_1,\delta_2,\delta_3,\delta_i) = \mbox{slope}([A]\delta_i)\in\bQ,$$
for any $i\geq4$. So, the normal lines $[A]\delta_1$, ... $[A]\delta_k$ meet (non-trivially) the lattice $\bZ^2$. Thus, for any representative $A\in GL(2,\bR)$ of $[A]$, the lattice $A^{-1}\bZ^2$ intersects non-trivially each of the lines $\delta_1$, ... $\delta_k$.
\end{proof}

\begin{corollary} \label{condBirapport1}
A quadrilateral is of rational type if and only if the cross ratio of its normals is rational or infinite. Moreover, the set of quadrilaterals of rational type is dense in the family of quadrilaterals and contains connected subfamilies.
\end{corollary}

\begin{corollary} \label{condBirapport2}
An orthotoric polytope with parameters $\beta_1<\beta_2<\alpha_1<\alpha_2$ is of rational type if and only $\cratio=\frac{(\beta_2-\alpha_1)(\alpha_2-\beta_1)}{(\beta_2- \beta_1)(\alpha_2-\alpha_1)}$ is rational.
\end{corollary}

Consider the orthotoric quadrilateral $\Delta_{\alpha,\beta} =\sigma([0,\beta]\times[1,\alpha])$, where $\sigma(x,y)=(x+y,xy)$, given by the characteristic pair $(\alpha,\beta)$ with $0<\beta<1<\alpha$, $\alpha-\beta\geq 1$, see Corollary~\ref{coroCharactPair}. The condition of being of rational type then read as $$\cratio(\alpha,\beta) = \frac{\alpha(\beta-1)}{\beta(\alpha-1)}\in\bQ.$$

\begin{proposition}\label{expressdesC}A labeled orthotoric quadrilateral associated to orthotoric parameters $(\alpha_1,\alpha_2,\beta_1,\beta_2,C_{\alpha_1}, C_{\alpha_2}, C_{\beta_1},C_{\beta_2})$ is a rational labeled polytope if and only if
\begin{itemize}
  \item[(1)] $\cratio=\frac{(\beta_2-\alpha_1)(\alpha_2-\beta_1)}{(\beta_2- \beta_1)(\alpha_2-\alpha_1)}\in \bQ,$
  \item[(2)] $C_{\beta_2}>0$ and there exist positive rational numbers $p_{\beta_1}$, $p_{\alpha_2}$, $p_{\alpha_1}$ such that $p_{\beta_1}C_{\beta_1} = \frac{(\beta_2-\alpha_1)}{(\alpha_1-\beta_1)}C_{\beta_2}$, $p_{\alpha_2}C_{\alpha_2} = -\frac{(\beta_2-\beta_1)}{(\alpha_2-\beta_1)}C_{\beta_2}$ and $p_{\alpha_1}C_{\alpha_1} = \frac{(\beta_2-\beta_1)}{(\alpha_1-\beta_1)}C_{\beta_2}$.
\end{itemize}
\end{proposition}
\begin{proof} In order to prove the proposition, we start with the following easy lemma.
\begin{lemma} \label{rseaucarct}
Let $u_0$, $u_1$, $u_2$ be pairwise linearly independent vectors of a $2$--dimensional vector space $V$. They generate a lattice if and only if there exist non-zero integers $n_0$, $n_1$, $n_2$ such that $n_0u_0 + n_1u_1 +n_2u_2=0$.
\end{lemma}

Using Lemma~\ref{rseaucarct} for both $\{u_{C_{\alpha_1}}, u_{C_{\alpha_2}}, u_{C_{\beta_1}}\}$ and $\{u_{C_{\alpha_1}}, u_{C_{\beta_2}}, u_{C_{\beta_1}}\}$, we obtain the homogeneous linear system
\begin{equation*}
  \begin{split}
    & n_1\alpha_1 C_{\alpha_1} + n_2\alpha_2C_{\alpha_2} + n_0\beta_1C_{\beta_1}=0 \;\;\; \;\;\;  n_1C_{\alpha_1} + n_2C_{\alpha_2} + n_0C_{\beta_1}=0\\
    & k_1\alpha_1 C_{\alpha_1} + k_2\beta_2C_{\beta_2} + k_0\beta_1C_{\beta_1}=0 \;\;\; \;\;\; k_1C_{\alpha_1} + k_2C_{\beta_2} + k_0C_{\beta_1}=0
  \end{split}
\end{equation*}
 for the unknowns $C_{\alpha_1}$, $C_{\alpha_2}$, $C_{\beta_1}$, $C_{\beta_2}$. It follows $C_{\beta_1} = \frac{k_2}{k_0}\frac{(\beta_2-\alpha_1)}{(\alpha_1-\beta_1)}C_{\beta_2}$,
$C_{\alpha_2} = -\frac{n_0k_2}{k_0n_2}\frac{(\beta_2-\alpha_1)}{(\alpha_1-\alpha_2)}C_{\beta_2}$, $C_{\alpha_1} = \frac{k_2}{k_1}\frac{(\beta_2-\beta_1)}{(\alpha_1-\beta_1)}C_{\beta_2}$ and
$$C_{\alpha_1} = \frac{k_2n_0}{k_0n_1}\frac{(\alpha_2-\beta_1)(\beta_2-\alpha_1)}{(\alpha_1-\beta_1)(\alpha_1-\alpha_2)}C_{\beta_2}= \cratio\frac{k_2n_0}{k_0n_1}\frac{(\beta_2-\beta_1)}{(\alpha_1-\beta_1)}C_{\beta_2},$$
from where we get the cross ratio condition (1) (since $\cratio=\frac{n_1k0}{k_1n_0}$). The expression of the coefficients of condition (2) follows easily.

Conversely, if conditions (1) and (2) are satisfied, the equations \begin{equation*}
    p_{\beta_1}u_{C_{\beta_1}} + p_{\alpha_1}u_{C_{\alpha_1}}  +  u_{C_{\beta_2}} =0\;\mbox{ and }\;
     p_{\beta_1}u_{C_{\beta_1}} + \cratio p_{\alpha_1}u_{C_{\alpha_1}} +  \cratio p_{\alpha_2}u_{C_{\alpha_2}} =0
\end{equation*} have rational coefficients. Then $u_{C_{\beta_1}}$, $u_{C_{\alpha_1}}$, $u_{C_{\beta_2}}$, $u_{C_{\alpha_2}}$ are all contained in a lattice.\end{proof}

From Proposition~\ref{crossratioTHEO} we know that any trapezoid is of rational type. However, $\Vlabel$s of a trapezoid must satisfy some condition  in order to be contained in a lattice. The following proposition gives these conditions. The proof is similar to the proof of Proposition~\ref{expressdesC}.
\begin{proposition}\label{expressdesCcalabi}
A labeled Calabi trapezoid with Calabi parameters $$(\alpha_1,\alpha_2,\beta_1,\beta_2, C_{\alpha_1}, C_{\alpha_2}, C_{\beta_1},C_{\beta_2})$$ is a rational labeled polytope if and only if $C_{\beta_2}>0$ and there exist positive rational numbers $p_{\beta_1}$, $p_{\alpha_2}$, $p_{\alpha_1}$ such that $-p_{\beta_1}C_{\beta_1} = C_{\beta_2}$, $p_{\alpha_2}\alpha_2C_{\alpha_2} = (\beta_2-\beta_1)C_{\beta_2}$ and $-p_{\alpha_1}\alpha_1C_{\alpha_1} = -(\beta_2-\beta_1)C_{\beta_2}$. \end{proposition}

 \begin{lemma} \label{lemstronglyRatIFFrat} $\Delta$ is strongly rational if and only if $\alpha$, $\beta \in\bQ$.
 \end{lemma}
\begin{proof} Recall that a polytope $\Delta$ sitting in a vector space $V$ is strongly rational if there exists a lattice $\Lambda^*\subset V$ such that, up to translation, all the vertices of $\Delta$ lie in the lattice $\Lambda^*$. Notice that, in this case, if one vertex lies in $\Lambda^*$ then all the vertices do and, seen as vectors with respect to the origin, the vertices generate a sublattice. Suppose that $\Delta$ is a quadrilateral and take a normal form of $\Delta$. The vertices $(0,0)$, $(0,1)$, $(1,0)$ and $(\alpha,1-\beta)$ belong to one lattice if and only if $\alpha$, $\beta \in\bQ$.\end{proof}

\subsection{Existence of extremal orthotoric and Calabi toric metrics} \label{sectOBSexistTHEO}

 Corollary~\ref{theoExistKE} from the introduction is a particular case of the following more general result. \begin{proposition}\label{theoExistFamily} Let $\Delta$ be a strongly rational convex quadrilateral which is not a parallelogram.
\begin{itemize}
   \item If $\Delta$ is generic, there exists a family, parameterized by $3$ positive rational numbers, of unstable symplectic toric orbifolds admitting no compatible extremal metric and whose moment polytope is $\Delta$.
   \item There exists a family, parameterized by $3$ positive rational numbers, of orthotoric extremal K\"ahler orbifolds whose moment polytope is $\Delta$. Moreover, this family contains a $2$--parameter subfamily of constant scalar curvature K\"ahler orbifolds and a $1$--parameter (sub-)subfamily of homothetic K\"ahler--Einstein orbifolds.
 \end{itemize} \end{proposition}

\begin{proof}Let $\Delta$ be a convex quadrilateral which is not a parallelogram. Denote by $\bfN(\Delta)$ the $4$--dimensional cone of inward normals associated to the facets of $\Delta$. We fix inward normal vectors $u= (u_1,u_2, u_3, u_4)$ so that any other normal inward vectors can be expressed as $u(r)$ for some $r\in\bR^4_{> 0}$, using conventions (\ref{setofLabel}) and~(\ref{NotationCr}). Let us define \begin{equation*}
 \bfR(\Delta) = \{r\in \bR_{>0}^4\,|\; (\Delta,u(r)) \mbox{ is a rational labeled polytope}\},
\end{equation*} so that the extremal orthotoric orbifolds with moment polytope $\Delta$ are in bijective correspondence with the elements of $\bfB^+(\Delta)\cap \bfR(\Delta)$. Unstable toric orbifolds with moment polytope $\Delta$ are in bijective correspondence with elements of $(\bfB(\Delta)\backslash\bfB^+(\Delta) )\cap \bfR(\Delta)$. Similarly, cscK (resp. KE) toric orbifolds with moment polytope $\Delta$ are in bijective correspondence with the points of $\bfS^+(\Delta)\cap \bfR(\Delta)$ (resp. $\bfK^+(\Delta)\cap \bfR(\Delta)$). Proposition~\ref{theoExistFamily} follows then from Lemma~\ref{Rat+KE}. \end{proof}

\begin{lemma}\label{Rat+KE} If $\Delta$ is strongly rational then $\bfR(\Delta)$ contains dense subsets of $\bfB^+(\Delta)$, $\bfC(\Delta)$, $\bfK(\Delta)$ and, if $\Delta$ is generic, $\bfB(\Delta)\backslash\bfB^+(\Delta)$.\end{lemma}

\begin{proof} Recall that a convex quadrilateral determines and is determined by its characteristic pair $(\alpha,\beta)\in\bR^2$ with $0\leq\beta<1<\alpha$ and $\alpha-\beta \geq 1$, see Corollary~\ref{coroCharactPair}. Thanks to Lemma~\ref{lemstronglyRatIFFrat} we know that under the hypothesis of the lemma (that is $\Delta$ is strongly rational) $\alpha$, $\beta \in\bQ$.

Suppose first that $\Delta$ is generic. Recall that $\bfR(\Delta)\neq \emptyset$ if and only if $\cratio(\alpha,\beta) = \frac{\alpha(\beta-1)}{\beta(\alpha-1)}\in\bQ$, see Corollary~\ref{condBirapport2}. Moreover, thanks to Proposition~\ref{expressdesC}, if $\bfR(\Delta) \neq \emptyset$ then
\begin{align*}\bfR(\Delta)=
\left\{\left(\frac{s q_1}{\beta}, \frac{\alpha s q_2}{\beta}, \frac{s q_3}{1-\beta},sq_4\right) \right.\left|\left.\begin{array}{l}  q_1,q_2,q_3,q_4 \in\bQ_{>0},\\
s\in\bR_{>0}\end{array}\right.\right\}.
\end{align*} If $\Delta$ is a trapezoid, then, thanks to Proposition~\ref{expressdesCcalabi}, $\bfR(\Delta) \neq \emptyset$ and
\begin{align*}\bfR(\Delta)=
\left\{\left( sq_1, s q_2, s q_3, sq_4\right) \right.\left|\left.\begin{array}{l}  q_1,q_2,q_3,q_4 \in\bQ_{>0},\\
s\in\bR_{>0} \end{array}\right.\right\}.
\end{align*}

Hence, Lemma~\ref{Rat+KE} follows from Theorem~\ref{propCONEsol} together with the fact that if $\alpha$, $\beta \in\bQ$ then $\bfR(\Delta)$ contains $\bQ_{>0}^4$ and the equations defining $\bfB(\Delta)$, $\bfC(\Delta)$ and $\bfK(\Delta)$ have rational coefficients. \end{proof}

\begin{rem} The strong rationality is necessary. For instance, suppose that $\Delta$ is a generic polytope of rational type. The equations~(\ref{condEXTpoly5}) and~(\ref{anciennement-ii}) defining $\bfC(\Delta)\cap\bfR(\Delta)$ in $\bfN(\Delta)$ may be turned into equations involving polynomials of one variable with rational coefficients using the fact that $\cratio=\frac{\alpha(\beta-1)}{\beta(\alpha-1)}\in \bQ$ and the parametrization of $\bfR(\Delta)$ by rational numbers. In particular, the existence of a point in $\bfC(\Delta)\cap\bfR(\Delta)$ implies that $\alpha$, $\beta$ are algebraic of degree at most $3$. Notice also that the condition $\cratio\in \bQ$ implies that $\alpha$, $\beta$ have the same algebraic degree. Similarly, if $\Delta$ is rational but not strongly rational then $\bfK(\Delta)\cap\bfR(\Delta)$ is empty. This fact can also be inferred from the general theory since it is well-known that the moment polytope of a K\"ahler-Einstein orbifold is strongly rational.
\end{rem}

\bibliographystyle{abbrv}

\end{document}